\documentclass[12pt]{article}
\newcommand{\dst} {\displaystyle}

\newcommand {\Star}{\mathop{\rm Star}}
\newcommand {\Spiral}{\mathop{\rm Spiral}}
\newcommand {\Stt}{\mathop{\rm Star}[\tau]}
\newcommand {\Spt}{\mathop{\rm Spiral}[\tau]}
\newcommand {\Stte}{\mathop{\rm Star}[\tau,\eta]}
\newcommand {\Spte}{\mathop{\rm Spiral}[\tau,\eta]}

\newcommand {\gtp}{\mathop{\mathcal{G}^+[\tau]}\nolimits}
\newcommand {\gtep}{\mathop{\mathcal{G}^+[\tau,\eta]}\nolimits}
\newcommand {\gup}{\mathop{\mathcal{G}^+[1]}\nolimits}

\newcommand {\Aut}{\mathop{\rm Aut}\nolimits}

\newcommand {\aut}{\mathop{\rm aut}\nolimits}

\newcommand {\Hol}{\mathop{\rm Hol}\nolimits}
\renewcommand {\Im}{\mathop{\rm Im}\nolimits}
\renewcommand {\Re}{\mathop{\rm Re}\nolimits}

\newcommand {\G}{{\cal G}}

\newcommand {\Ss}{\mathcal{S}}
\newcommand {\Ll}{\mathcal{L}}

\newfont{\bbb}{msbm10 at 12pt}
\def\Bbb#1{\hbox{\bbb #1}}

\newcommand{\C}{{\Bbb C}}
\newcommand{\R}{{\Bbb R}}
\newcommand{\N}{{\Bbb N}}

\usepackage{amsfonts}
\usepackage{latexsym}
\usepackage{epsf,graphicx,latexsym,%
}

\newcommand{\pl}{\partial}
\newcommand{\pr}{\noindent{\bf Proof.}\quad }
\newcommand{\epr}{\ $\Box$ \vspace{2mm} }

\newcommand{\be} {\begin{eqnarray}}
\newcommand{\ee} {\end{eqnarray}}
\newcommand{\bep} {\begin{eqnarray*}}
\newcommand{\eep} {\end{eqnarray*}}

\newtheorem{defin}{Definition}
\newtheorem{theorem}{Theorem}
\newtheorem{corol}{Corollary}
\newtheorem{lemma}{Lemma}

\begin{document}

\title{A flower structure of backward flow invariant domains for
semigroups%
\thanks{\footnotesize {\it 2000 Mathematics Subject Classification:} 37C10, 30C45}}%

\author{Mark Elin \\ {\small Department of Mathematics, ORT  Braude College,} \\
{\small P.O. Box 78, Karmiel 21982, ISRAEL} \\ {\small e-mail:
mark.elin@gmail.com}
\\ David Shoikhet \\ {\small Department of Mathematics, ORT Braude College,}
\\ {\small P.O. Box 78, Karmiel 21982, ISRAEL} \\ {\small e-mail: davs27@netvision.net.il}
\\ Lawrence Zalcman\footnote{Research supported by The German--Israeli
Foundation for Scientific Research and Development, G.I.F. Grants
No. G-643-117.6/1999 and I-809-234.6/2003}
\\ {\small Department of Mathematics, Bar-Ilan University,} \\
{\small 52900 Ramat-Gan, ISRAEL} \\ {\small e-mail:
zalcman@macs.biu.ac.il}  }
\date{}


\maketitle

\begin{abstract}
In this paper, we study conditions which ensure the existence of
backward flow invariant domains for semigroups of holomorphic
self-mappings of a simply connected domain $D$. More precisely,
the problem is the following. Given a one-parameter semigroup
$\mathcal S$ on $D$, find a simply connected subset $\Omega\subset
D$ such that each element of $\mathcal S$ is an automorphism of
$\Omega$, in other words, such that $\mathcal S$ forms a
one-parameter group on $\Omega$.

On the way to solving this problem, we prove an angle distortion
theorem for starlike and spirallike functions with respect to
interior and boundary points.
\end{abstract}


Let $D$ be a simply connected domain in the complex plane $\C$. By
$\Hol(D,\Omega)$ we denote the set of all holomorphic functions on
$D$ with values in a domain $\Omega $ in $\C$. We write $\Hol(D)$
for $\Hol(D,D)$, the set of holomorphic self-mappings of $D$. This
set is a topological semigroup with respect to composition. We
denote by $\Aut(D)$ the group of all automorphisms of $D$; thus
$F\in\Aut(D)$ if and only if $F$ is univalent on $D$ and $F(D)=D.$

\begin{defin}
A family $\Ss=\left\{F_t\right\}_{t\geq 0}\subset\Hol(D)$ is said
to be a one-parameter continuous semigroup (semiflow) on $D$ if

(i) $F_{t}(F_{s}(z))=F_{t+s}(z)$ for all $t,s\geq 0,$

(ii) $\lim\limits_{t\rightarrow 0^+}F_{t}(z)=z$ for all $z\in D$.

\noindent If, in addition, condition (i) holds for all $t,s\in\R$,
then $(F_{t})^{-1}=F_{-t}$ for each $t\in\R$; and $\Ss$ is called
a {\it one parameter continuous group (flow)} on $D$. In this
case, $\Ss\subset\Aut(D)$.
\end{defin}

In this paper, we study the following problem. {
\bf Given a one-parameter semigroup $\Ss\subset\Hol(D)$, find a
simply connected domain $\Omega\subset D$ (if it exists) such that
$\Ss\subset\Aut(\Omega).$}

It is well-known that condition (ii) and holomorphy, in fact,
imply that
\[
\lim_{t\rightarrow s}F_{t}(z)=F_{s}(z)
\]
for each $z\in D$ and $s>0$ ($s\in\R$ in the case when
$\Ss\subset\Aut(D)$); see, for example, \cite{B-P}, \cite{AM-92},
\cite{R-S-96} and \cite{R-S-97}. This explains the name
``continuous semigroup" in our terminology.

Furthermore, it follows by a result of E. Berkson and H. Porta
\cite{B-P} that each continuous semigroup is differentiable in
$t\in\R^{+}=[0,\infty ),$ (see also \cite{AM-88} and
\cite{R-S-98}). So, for each continuous semigroup (semiflow)
$\Ss=\left\{F_{t}\right\}_{t\geq 0}\subset\Hol(D)$, the limit
\begin{equation}
\lim_{t\rightarrow 0^{+}}\frac{z-F_{t}(z)}{t}=f(z),\quad z\in D,
\label{1a}
\end{equation}
exists and defines a holomorphic mapping $f\in\Hol(D,\C)$. This
mapping $f$ is called the { \it (infinitesimal) generator of}
$\Ss=\left\{F_{t}\right\}_{t\geq 0}.$ Moreover, the function
$u(=u(t,z)),\ (t,z)\in\R^{+}\times D$, defined by
$u(t,z)=F_{t}(z)$ is the unique solution of the Cauchy problem
\be\label{2a}
\left\{
\begin{array}{l}
{\displaystyle\frac{\pl u(t,z)}{\pl t}}+f(u(t,z))=0,\vspace{3mm}\\
u(0,z)=z,\quad z\in D.
\end{array}
\right.
\ee

Conversely, a mapping $f\in\Hol(D,\C)$ is said to be a {\it
semi-complete (respectively, complete) vector field} on $D$ if the
Cauchy problem (\ref{2a}) has a solution $u(=u(t,z))\in D$ \ for
all $z\in D$ and $t\in\R^{+}$ (respectively, $t\in\R$). Thus
$f\in\Hol(D,\C)$ is a semi-complete vector field if and only if it
is the generator of a one-parameter continuous semigroup $\Ss$
(semiflow) on $D$. It is complete if and only if
$\Ss\subset\Aut(D).$ The set of semi-complete vector fields on $D$
is denoted by $\G(D)$. The set of complete vector fields on $D$ is
usually denoted by $\aut(D)$ (see, for example, \cite{I-S},
\cite{UH}, \cite{SD}).

Thus, in these terms, our problem can be rephrased as follows.
{\bf Given $f\in\G(D)$, find a domain $\Omega$ (if it exists) such
that $f\in\aut(\Omega)$}.

Let now $D=\Delta $ be the open unit disk in $\C$. In this case,
$\G(\Delta)$ is a real cone in $\Hol(\Delta,\C)$, while
$\aut(\Delta)\subset\G(\Delta)$ is a real Banach space (see, for
example, \cite{R-S-98}). Moreover, by the Berkson--Porta
representation formula, a function $f$ belongs to $\G(\Delta)$ if
and only if there is a point $\tau\in\overline{\Delta}$ and a
function $p\in\Hol(\Delta,\C)$ with positive real part ($\Re
p(z)\ge0$ everywhere) such that
\begin{equation}
f(z)=(z-\tau )(1-z\overline{\tau })p(z).  \label{3a}
\end{equation}
This representation is unique and is equivalent to
\[
f(z)=a-\bar{a}z^2+zq(z),\quad a\in\C,\ \Re q(z)\ge0
\]
(see \cite{A-E-R-S}). Moreover, $f\in\Hol(\Delta,\C)$ is complete
if and only if it admits the representation
\begin{equation}
f(z)=a-\bar{a}z^{2}+ibz  \label{4a}
\end{equation}
for some $a\in\C$ and $b\in\R$ (see, \cite{B-K-P}, \cite{AJ},
\cite{UH}).

Note also that if a semigroup $\Ss=\left\{F_t\right\}_{t\geq 0}$
generated by $f\in\G(\Delta)$ does not contain an elliptic
automorphism of $\Delta$, then the point
$\tau\in\overline{\Delta}$ in representation (\ref{3a}) is the
unique attractive point for the semigroup $\Ss$, i.e.,
\begin{equation}\label{5a}
\lim_{t\to\infty}F_t(z)=\tau
\end{equation}
for all $z\in \Delta$. This point is usually referred as the
Denjoy--Wolff point of $\Ss$. In addition,

$\bullet$ if $\tau\in\Delta,$ then $\tau =F_{t}(\tau )$ is a
unique fixed point of $\Ss$ in $\Delta$;

$\bullet$ if $\tau\in \partial \Delta ,$ then
\[
\tau =\lim_{r\rightarrow 1^{-}}F_{t}(r\tau )
\]
is a common boundary fixed point of $\Ss$ in $\overline{\Delta}$,
and no element $F_t\ (t>0)$ has an interior fixed point in
$\Delta.$

Also, we observe that for $\tau\in\Delta$, formula (\ref{3a})
implies the condition
\begin{equation}
\Re f'(\tau )\geq 0.  \label{6a}
\end{equation}
Comparing this with (\ref{3a}) and (\ref{4a}), we see that $\Ss$
consists of elliptic automorphisms if and only if
\begin{equation}
\Re f'(\tau )=0.  \label{7a}
\end{equation}
Consequently, condition (\ref{5a}) is equivalent to
\begin{equation}\label{8a}
\Re f'(\tau )>0.
\end{equation}

If $\tau$ in (\ref{3a}) belongs to $\pl \Delta,$ then if follows
by the Riesz--Herglotz representation of the function $p$ in
(\ref{3a}) that the angular limits
\begin{equation}\label{9a}
f(\tau):=\angle\lim_{z\to\tau}f(z)=0\quad\mbox{and}\quad
f'(\tau):=\angle\lim_{z\rightarrow \tau }f'(z)=\beta
\end{equation}
exist and that $\beta$ is a nonnegative real number (see also
\cite{E-S1}). Moreover, if for some point $\zeta\in\pl\Delta$
there are limits
\[
\angle \lim_{z\to\zeta}f(z)=0
\]
and
\[
\angle \lim_{z\to\zeta}f'(z)=\gamma
\]
with $\gamma\geq0,$ then $\gamma=\beta$ and $\zeta=\tau$ (see
\cite{E-S1} and \cite{SD1}).

In the case where $\beta >0$, the semigroup
$\Ss=\left\{F_t\right\}_{t\geq 0}$ consists of mappings
$F_{t}\in\Hol(\Delta)$ of hyperbolic type,
\[
\angle \lim_{z\to\tau}\frac{\pl F_t(z)}{\pl z} =e^{-t\beta }<1;
\]
otherwise $(\beta =0)$, it consists of mappings of parabolic type,
\[
\angle \lim_{z\to\tau}\frac{\pl F_t(z)}{\pl z} =1\quad \mbox{for
all }t\geq 0.
\]

For $\tau\in\overline\Delta$, we use the notation $\gtp$ for a
subcone of $\G(\Delta)$ of functions $f$ defined by (\ref{3a}) for
which
\begin{equation}
\Re f'(\tau )>0.  \label{10a}
\end{equation}
We solve the problem mentioned above for the class $\gtp$ of
generators.


\begin{defin}
Let $\Ss=\{F_{t}\}_{t\geq 0}$ be a semiflow on $\Delta$. A domain
$\Omega\subset\Delta$ is called a (backward) {\bf flow-invariant
domain} (shortly, {\bf FID}) for $\Ss$ if
${\Ss\subset\Aut(\Omega)}.$
\end{defin}

We need the following notation. We write $f\in\gtep$, where
${\tau\in\overline\Delta}$, ${\eta\in\pl\Delta,}\ \eta\not=\tau$,
if $f\in\gtp$, $f(\eta)=\angle \lim\limits_{z\to\eta}f(z)=0$ and
$\gamma=\angle\lim\limits_{z\to\eta}f'(z)$ exists finitely. In
fact, in this case $\gamma$ must be a real negative number (see
Lemma~6 below).

\begin{theorem}
Let $\Ss=\{F_{t}\}_{t\geq 0}$ be a semiflow on $\Delta$ generated
by $f\in\gtp$, for some $\tau\in\overline\Delta$ with $f(\tau)=0$
and $f'(\tau)=\beta, \ \Re\beta>0$. The following assertions are
equivalent.

(i) $f\in\gtep$ for some $\eta\in\pl\Delta$.

(ii) There is a nonempty (backward) flow invariant domain
$\Omega\subset\Delta$, so $\Ss\subset\Aut(\Omega)$.

(iii) For some $\alpha>0$, the differential equation
\begin{equation}\label{f**}
\alpha \varphi'(z)(z^{2}-1)=2f(\varphi(z))
\end{equation}
has a locally univalent solution $\varphi$ with $|\varphi(z)|<1$
when $z\in\Delta$. Moreover, in this case $\varphi$ is univalent
and is a Riemann mapping of $\Delta$ onto a flow invariant domain
$\Omega$.
\end{theorem}
This theorem can be completed by the following result.
\begin{theorem}
Let $\Ss=\{F_{t}\}_{t\geq 0}$ be a semiflow on $\Delta$ generated
by $f\in\gtp$, for some $\tau\in\overline\Delta$ with $f(\tau)=0$
and $f'(\tau)=\beta, \ \Re\beta>0$. The following assertions hold.

(a) If $f\in\gtep$ for some $\eta\in\pl\Delta$ with
$\gamma=\angle\lim\limits_{z\to\eta}f'(z)$, then for each
$\alpha\ge-\gamma$, equation {\rm (\ref{f**})} has a univalent
solution $\varphi$ such that $\varphi(1)=\tau,\ \varphi(-1)=\eta$
and $\Omega=\varphi(\Delta)$ is a (backward) flow invariant domain
for $\Ss$. In addition, $\tau
=\lim\limits_{t\to\infty}F_t(z)\in\pl\Omega,\ z\in\Omega$, and
$\lim\limits_{t\to-\infty}F_t(z)=\eta\in\pl\Delta\cap\pl\Omega$
for each $z\in\Omega$.

(b) If $\Omega\subset\Delta$ is a nonempty (backward) flow
invariant domain, then it is a Jordan domain such that
$\tau\in\pl\Omega$, and there is a point
$\eta\in\pl\Omega\cap\pl\Delta$ such that
$\lim\limits_{t\to-\infty}F_t(z)=\eta$ whenever $z\in\Omega$,
$\angle\lim\limits_{z\to\eta}f(z)=0$ and
$\angle\lim\limits_{z\to\eta}f'(z)=:\gamma$ exists with
$\gamma<0$. In addition, there is a conformal mapping $\varphi$ of
$\Delta$ onto $\Omega$ which satisfies equation {\rm (\ref{f**})}
with some $\alpha\ge-\gamma$.

(c) Conversely, if for some $\alpha>0$, the differential equation
{\rm(\ref{f**})} has a locally univalent solution
$\varphi\in\Hol(\Delta)$, then it is, in fact, a conformal mapping
of $\Delta$ onto the FID \quad ${\Omega=\varphi(\Delta)}$ such
that $\varphi(1)=\tau\in\pl\Omega$ and $\varphi(-1)=\eta$ for some
$\eta\in\pl\Delta\cap\pl\Omega$.

In addition, $f(\eta)=0$ and $f'(\eta)=\gamma$ with
$0>\gamma\ge-\alpha$.
\end{theorem}

\begin{defin}
A (backward) flow-invariant domain (FID) $\Omega\subset\Delta$ for
$\Ss$ is said to be maximal if there is no
$\Omega_1\supset\Omega,\ \Omega_1\neq \Omega ,$ such that
$\Ss\subset\Aut(\Omega _{1}).$
\end{defin}

\begin{theorem}
Let $f\in\gtep$ for some $\tau\in\overline\Delta,\
\eta\in\pl\Delta$ with $\gamma=f'(\eta)\Bigl(<0\Bigr)$, and let
$\varphi$ be a (univalent) solution of {\rm(\ref{f**})} with some
$\alpha\ge-\gamma$ normalized by $\varphi(1)=\tau$ and
$\varphi(-1)=\eta$. The following assertions are equivalent:

(i) $\Omega=\varphi(\Delta)$ is a maximal FID;

(ii) $\alpha=-\gamma$;

(iii) $\varphi$ is isogonal at the boundary point $z=-1$ (see
Remark 3 below).
\end{theorem}

\noindent{\bf Remark 1.} In general, a maximal FID for $\Ss$ need
not be unique. Theorem~1 states that if $\Ss=\{F_{t}\}_{t\geq 0}$
is generated by $f\in\gtp$, then its {\it FID is not empty if and
only if there is a point $\eta\in\pl\Delta,$ such that
$f(\eta)=\angle \lim\limits_{z\to\eta}f(z)=0$ and
$f'(\eta)=\angle\lim\limits_{z\to\eta}f'(z)$ exists finitely with
$f'(\eta)<0$.} This point $\eta $ is a repelling fixed point for
$\Ss=\{F_{t}\}_{t\geq 0}$ as $t\to\infty$, namely,
$F_t(\eta)=\eta$ and $\left.\frac{\pl F_t(z)}{\pl
z}\right|_{z=\eta}=e^{-tf'(\eta)}>1$ (see \cite{E-S1}). Moreover,
there is {\it a one-to-one correspondence between maximal flow
invariant domains for $\Ss$ and such repelling fixed points.}

\begin{theorem}
Let $f\in\G^+[\tau,\eta_k]$ for some sequence
$\{\eta_k\}\in\pl\Delta$, i.e., $f(\eta_k)=0$ and
$\gamma_k=f'(\eta_k)>-\infty$.

The following assertions hold.

(i) There is $\delta>0$ such that $\gamma_k<-\delta<0$ for all
$k=1,2,\ldots$.

(ii) For each $a<-\delta<0$ there is at most a finite number of
the points $\eta_k$ such that $a\le\gamma_k<-\delta$.

Consequently equation (\ref{f**}) has a (univalent) solution
$\varphi\in\Hol(\Delta)$ for each
$\alpha\ge-\max\{\gamma_k\}>-\delta$.

(iii) If $\varphi_k$ is a solution of (\ref{f**}) normalized by
$\varphi_k(1)=\tau,\ \varphi_k(-1)=\eta_k$ with $\alpha=\gamma_k$
and $\Omega_k=\varphi_k(\Delta)$ (i.e., $\Omega_k$ are maximal),
then for each pair $\Omega_{k_1}$ and $\Omega_{k_2}$ such that
$\eta_{k_1}\neq\eta_{k_2}$ either
$\overline{\Omega_{k_1}}\cap\overline{\Omega_{k_2}}=\{\tau\}$ or
$\overline{\Omega_{k_1}}\cap\overline{\Omega_{k_2}}=l$, where $l$
is a continuous curve joining $\tau$ with a point on $\pl\Delta$.
\end{theorem}

We illustrate the content of our theorems in the following
examples.

\vspace{3mm}

\noindent{\bf Example 1.} Consider a generator $f\in\G^+[0]$
defined by
\[
f(z)=z(1-z^n),\quad n\in\N.
\]
Solving the Cauchy problem (\ref{2}), we find
\[
F_t(z)=\frac{ze^{-t}}{\sqrt[n]{1-z^n+z^ne^{-nt}}}\,.
\]
In this case, $f$ has $n$ additional null points
$\dst\eta_k=e^{\frac{2\pi ik}n},\,k=1,2,\ldots,n,$ on the unit
circle with finite angular derivative $\gamma=f'(\eta_k)=-n$. So
the generated semiflow has $n$ repelling fixed points, and there
are $n$ maximal flow invariant domains. One can show that the
functions
\[
\varphi_k(z)=e^{\frac{2\pi ik}n}\,\sqrt[n]{\frac{1-z}2}
\]
are the solutions of (\ref{f**}) with $\alpha=n$ satisfying
$\varphi_k(1)=0$ and $\varphi_k(-1)=\eta_k$ which map $\Delta$
onto $n$ FID's $\Omega_k$ (for $n=2$, these domains form
lemniscate) with
$\overline{\Omega_i}\bigcap\overline{\Omega_j}=\{0\}$ when $i\neq
j$. The family $\{F_t\}_{t\in\R}$ forms a group of automorphisms
of each one of these domains. See Figure~1 for $n=1,2,3$ and $5$.
For $n=1$, for instance, it can be seen explicitly that
$F_t(\varphi(z))$ is well-defined for all $t\in\R$ and tends to
$\eta=1$ when $t\to-\infty$.
\begin{figure}\centering 
    \includegraphics[angle=270,width=5.7cm,totalheight=5.7cm]{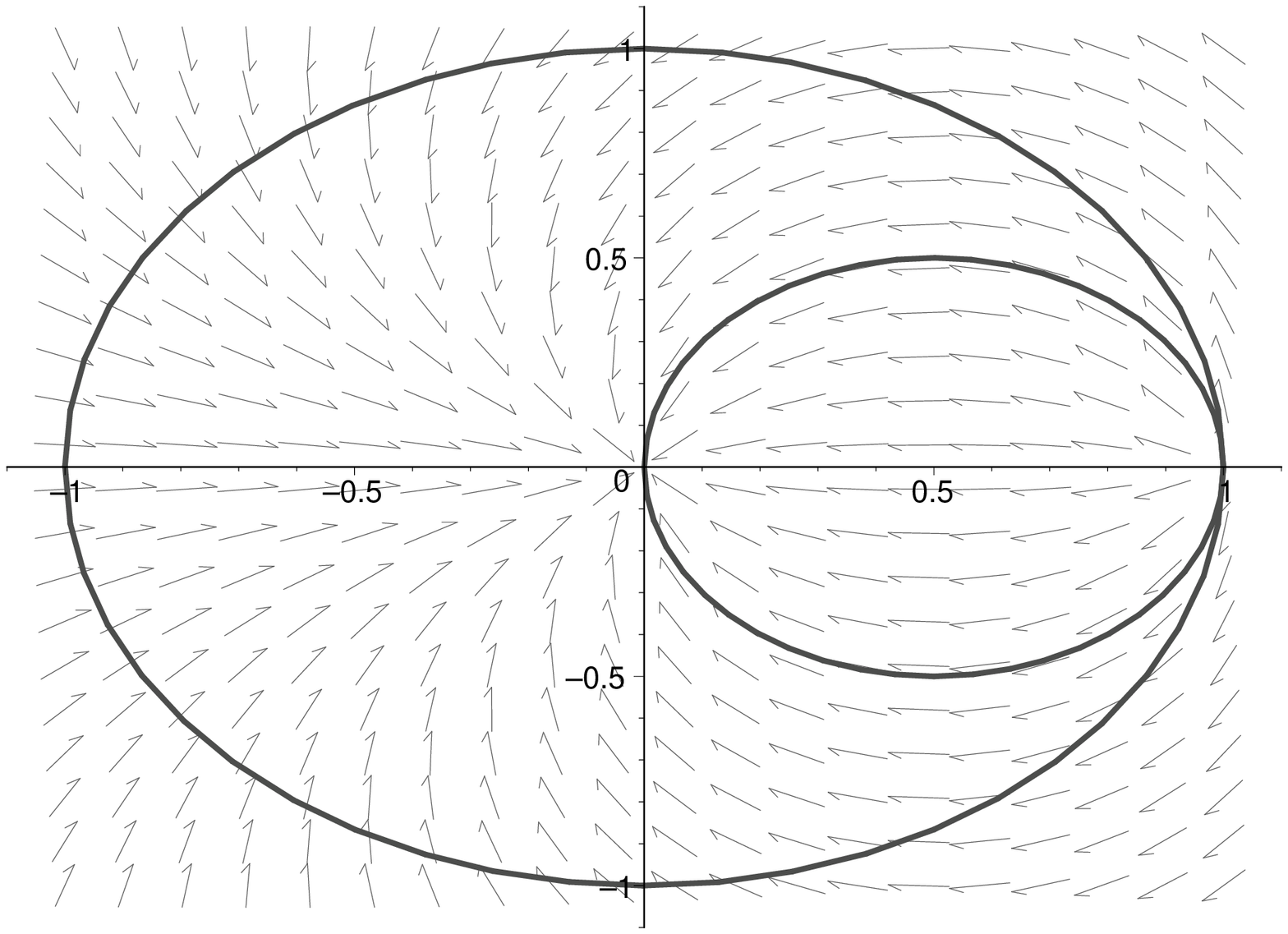}
    \hspace{7mm}
    \includegraphics[angle=270,width=5.7cm,totalheight=5.7cm]{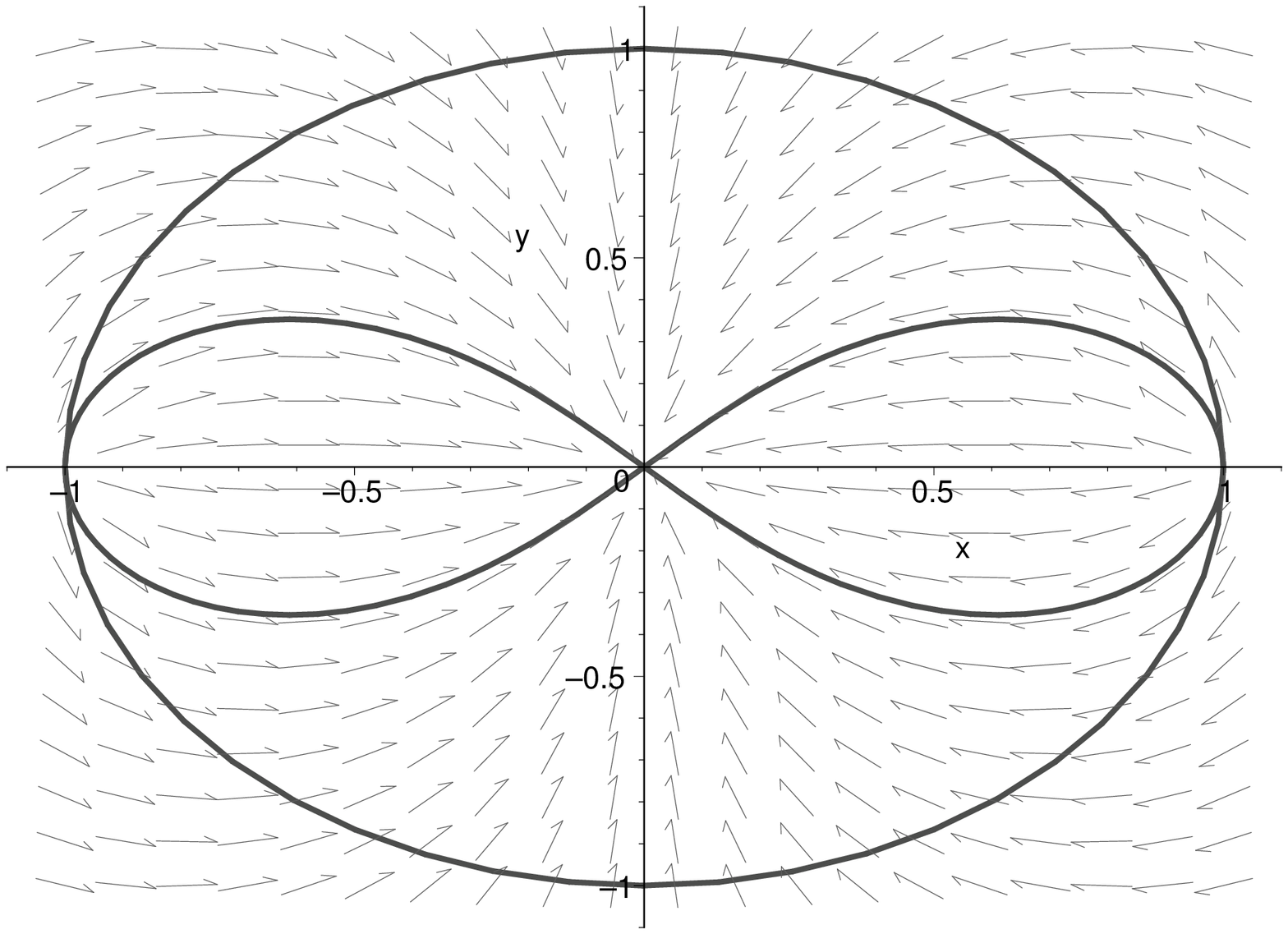}
\\ \vspace{5mm}
    \includegraphics[angle=270,width=5.7cm,totalheight=5.75cm]{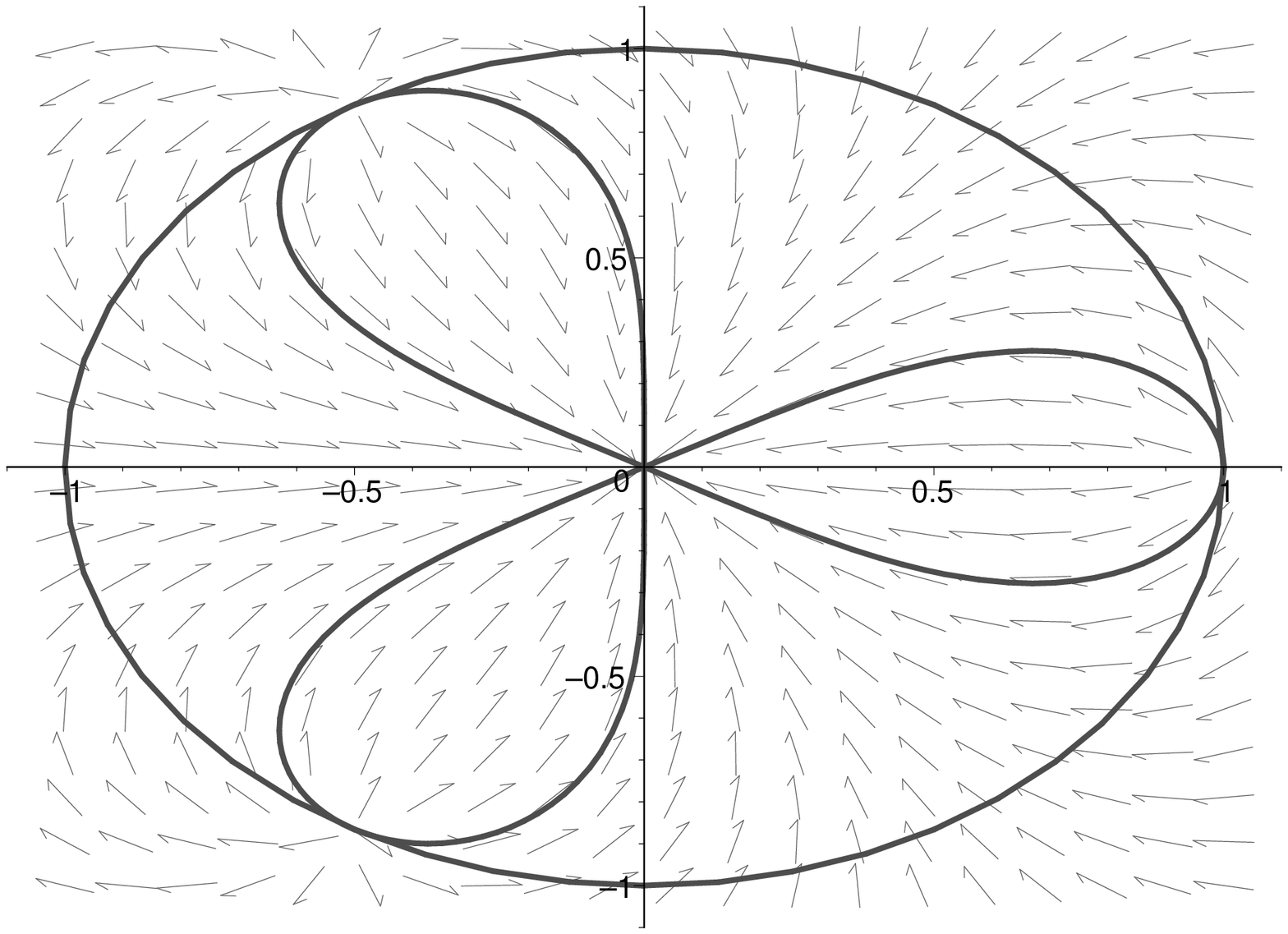}
\hspace{7mm}
    \includegraphics[angle=270,width=5.7cm,totalheight=5.7cm]{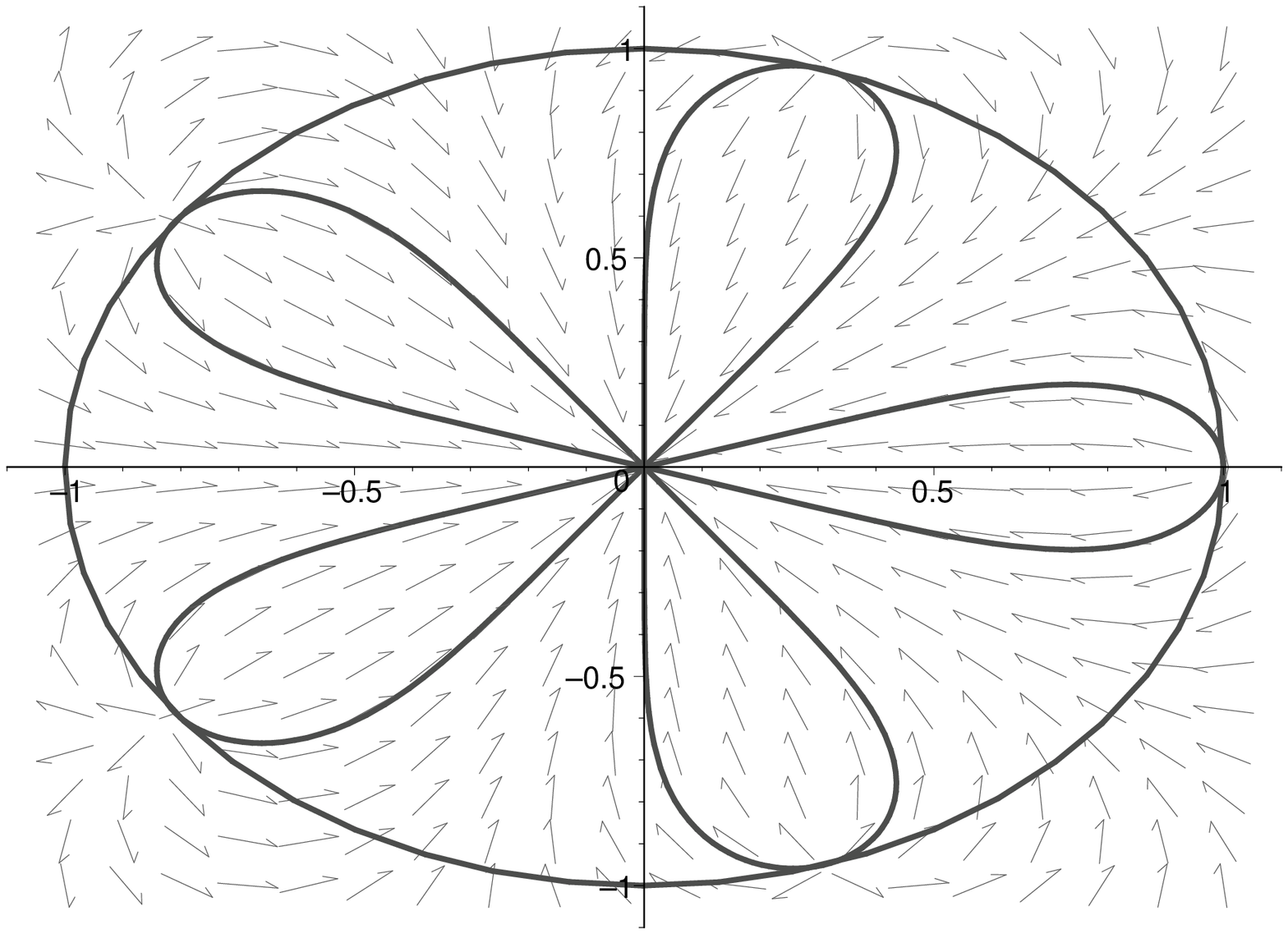}
    \caption{Example 1, $n=1,2,3,5$.}
\end{figure}

\vspace{2mm}

\noindent{\bf Example 2.} Consider a generator $f\in\G^+[1]$
defined by
\[
f(z)=-\frac{(1-z)(1+z^2)}{1+z}\,.
\]
Solving the Cauchy problem (\ref{2}), we find
\[
F_t(z)=\frac{(1+z^2)e^{2t}-(1-z)\sqrt{2(1+z^2)e^{2t}-(1-z)^2}}{(1+z^2)e^{2t}-(1-z)^2}\,.
\]
Since $f$ has the two additional null points $\eta_{1,2}=\pm
i\in\pl\Delta$ with finite angular derivative $\gamma=f'(\pm
i)=-2$, the generated semiflow has two repelling fixed points.
Thus, there are two maximal flow invariant domains $\Omega_1$ and
$\Omega_2$. One can show that these domains $\Omega_j$ coincide
with the upper and the lower half-disks (see Figure 2). So we have
$\overline{\Omega_1}\bigcap\overline{\Omega_2}=\{-1<x<1\}$. In
each of these two domains, the family $\{F_t\}_{t\in\R}$ is well
defined and forms a group of automorphisms.
\begin{center}
\begin{figure}[h]
\centering
\includegraphics[angle=270,width=6cm,totalheight=6cm]{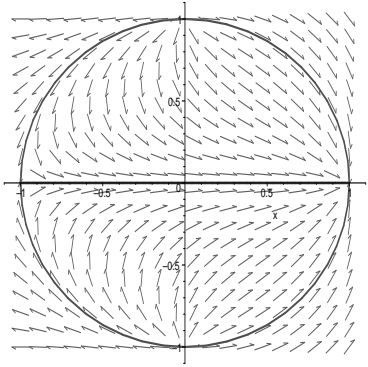}
\caption{Example 2. The flow generated by
$f(z)=-\frac{(1-z)(1+z^2)}{1+z}$ and two flow-invariant domains}
\end{figure}
\end{center}

\vspace{-5mm}

The following example shows that a maximal flow invariant domain
may be even dense in the open unit disk.

\noindent{\bf Example 3.} Let $f\in\G^+[0]$ be given by
\[
f(z)=z\,\frac{1-z}{1+z}\,.
\]
In this case, $\tau =0$ and $\eta =1$. Also, we have $f'(0)=1$ and
${f'(1)=-\frac{1}{2}\,.}$ Solving equation (\ref{f**}) with
$\alpha=\frac12$, one can write its solution in the form $\varphi
(z)=h^{-1}(h_{0}(z))$, where $h$ is the Koebe function $\dst
h(z)=\frac{z}{(1-z)^{2}}$ and
${h_{0}(z)=\dst\left(\frac{1-z}{1+z}\right)^2}$. We shall see
below that each solution of (\ref{f**}) has a similar
representation.

Thus $\varphi$ maps $\Delta$ onto the maximal flow invariant
domain $\Omega=\varphi(\Delta)=\Delta\setminus\{-1\leq x\leq 0\}$;
see Figure 3. (All the pictures were obtained by using the vector
field drawing tool in Maple~9.)
\begin{center}
\begin{figure}[h]
\centering
\includegraphics[angle=270,width=6cm,totalheight=6cm]{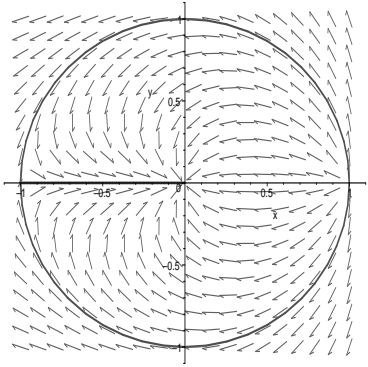}
\caption{Example 3. The flow generated by
$f(z)=z\frac{(1-z)}{1+z}$ and the dense flow-invariant domain}
\end{figure}
\end{center}

\noindent{\bf Remark 2.} Let $F\in\Hol(\Delta)$ be a single
self-mapping of $\Delta$ which can be embedded into a continuous
semigroup, i.e., there is a semiflow $\Ss=\{F_t\}_{t\geq 0}$ such
that $F=F_1$. In this case, all the fractional iterations $F_t$ of
$F$ have the same collection of boundary fixed points for all
$t\ge0$ (see \cite{C-DM-P}). In turn, our theorem asserts the
existence of backward fractional iterations of $F$ defined on a
FID $\Omega$ whenever $F$ has a repelling boundary fixed point
$\eta$, i.e.,
\be\label{rem1}
A=F'(\eta)=\lim_{z\to\eta}F'(z)>1.
\ee

As a matter of fact, for a single mapping which is not necessarily
embedded into a semiflow (not even necessarily univalent on
$\Delta$), the existence of backward integer iterations under
condition (\ref{rem1}) was proved in \cite{PP-C}. This fact has
provided the existence of conjugations near repelling points. More
precisely, the main result in \cite{PP-C} asserts that if $\eta
=1,\ \dst a=\frac{A-1}{A+1}$ and $\dst G(z)=\frac{z-a}{1-az}\,$,
then there is $\varphi\in\Hol(\Delta)$ with $\varphi(1)=1$ which
is a conjugation for $F$ and $G$, i.e.,
\[
\varphi (G(z))=F(\varphi (z)).
\]
However, for the case in which $F$ can be embedded into a
continuous semigroup $\Ss=\{F_{t}\}$, it is not clear whether
$\varphi$ is a conjugation for the whole semiflow $\Ss$ and the
flow produced by $G$.

It is natural to expect a more precise result under stronger
requirements. A direct consequence of the proof of our Theorem~1
is the following assertion for conjugations.

\begin{corol}
Let $F\subset\Hol(\Delta)$ be embedded into a semiflow
$\Ss=\{F_{t}\}_{t\geq 0}$ of hyperbolic type and let
$\eta\in\pl\Delta$ be a repelling fixed point of $F$ with
$A=F'(\eta)>1$.

Then for each $B\geq A$ and the automorphism
$G(=G_{B})\in\Aut(\Delta)$ defined by
\[
G(z)=\frac{z+b}{1+zb}\,,
\]
where $b=\frac{B-1}{B+1}\,$, there is a homeomorphism
$\varphi(=\varphi_B)$ of $\overline{\Delta},\
\varphi\in\Hol(\Delta)$, such that $\varphi(\eta)=-1$ and
\[
\varphi(G(z))=F(\varphi(z)),\quad z\in \Delta.
\]

Moreover, for all $t\in\R$ and $w\in\varphi(\Delta)$, the flow
$\{F_t(w)\}_{t\in\R}$ is well-defined with $F_1=F$ and
\[
F_t(\varphi(z))=\varphi(G_{t}(z)),\quad\mbox{for all\ }\ t\in\R,
\]
where
\[
G_t(z)=\frac{z+1+e^{-\alpha t}(z-1)}{z+1- e^{-\alpha
t}(z-1)},\quad t\in\R,
\]
with $\alpha =\log B.$

In addition, $\varphi_B(\Delta)\subseteq \varphi_A(\Delta)$, with
$\varphi_A(\Delta)=\varphi_B(\Delta)$ if and only if $A=B.$
\end{corol}

Our approach to construct conjugations is different from that used
in \cite{PP-C}.

The main tool of the proof of our theorems is a linearization
method for semigroups which uses the classes of starlike and
spirallike functions on $\Delta$.

\begin{defin}\label{defstar}
A univalent function $h$ is called spirallike (respectively,
starlike) on $\Delta$ if for some $\mu\in\C$ with $\Re\mu>0$
(respectively, $\mu\in\R$ with $\mu>0$) and for each point
$z\in\Delta$,
\be\label{defspiral}
\left\{ e^{-\mu t}h(z),\ t\geq 0\right\}\subset h(\Delta).
\ee
\end{defin}

In this case, we say that $h$ is $\mu$-spirallike.

Obviously, $0\in\overline{ h(\Delta)}$.

\noindent$\bullet$ If $0\in h(\Delta)$, (i.e., if there is a point
$\tau\in\Delta$ such that $h(\tau)=0)$, then $h$ is called {\it
spirallike (respectively, starlike) with respect to an interior
point.}

\noindent$\bullet$ If $0\not\in h(\Delta)$ (and hence $0\in\pl
h(\Delta)$), $h$ is called {\it spirallike (respectively,
starlike) with respect to a boundary point.} In this case, there
is a boundary point $\tau\in\pl\Delta$ such that
$h(\tau):=\angle\lim\limits_{z\to \tau}h(z)=0$ (see, for example,
\cite{E-R-S1}).

The class of spirallike (starlike) functions satisfying
$h(\tau)=0,\ \tau\in\overline\Delta,$ is denoted by $\Spt$
(respectively, $\Stt$).

It follows from Definition \ref{defstar} that a family
$\Ss=\{F_t\}_{t\ge0}$ of holomorphic self-mappings of the open
unit disk $\Delta$ defined by
\[
F_t(z):=h^{-1}\left(e^{-\mu t}h(z)\right)
\]
forms a semiflow on $\Delta$. Differentiating this semiflow at
$t=0^+$, one sees that $h$ is a solution of the differential
equation
\be\label{stareq}
\mu h(z)=h'(z)f(z),
\ee
where $f\in\gtp$ is the generator of $\Ss$. As a matter of fact,
the converse assertion also holds \cite{E-R-S1}, \cite{E-R-S2},
\cite{A-E-S}, \cite{E-Go-R-S}, \cite{E-R-S3}. More precisely, we
have

\begin{lemma}\label{1}
Let $\Ss=\{F_t\}_{t\ge0}$ be a semigroup of holomorphic
self-mappings generated by $f\in\gtp,\ \tau\in\overline\Delta$.

(i) If $\tau\in\Delta$, then equation (\ref{stareq}) has a
univalent solution if and only if $\mu=f'(\tau)$.

(ii) If $\tau\in\pl\Delta$, then equation (\ref{stareq}) has a
univalent solution $h$ satisfying $h(\tau)=0$ if and only if
$\mu\in\Lambda_\beta:=\left\{w\not=0:\ |w-\beta|\le\beta\right\}$,
where $\beta=f'(\tau)$.

Moreover, in both cases, this solution $h$ is a spirallike
(starlike) function which satisfies Schr\"oder's functional
equation
\be\label{schroder}
h(F_t(z))=e^{-\mu t}h(z),\quad t\ge0,\ z\in\Delta.
\ee
\end{lemma}

It is clear that $h$ is $\lambda$-spirallike for each $\lambda$
with $\arg\lambda=\arg\mu\in\left(-\frac\pi2, \frac\pi2\right)$.
We call this function $h$ the {\it spirallike (starlike) function
associated with $f$.}

Since we are interested in generators having additional null
points on the boundary, we introduce the following subclasses of
$\gtp$ and of $\Spt$ ($\Stt$).

\noindent$\bullet$ Given $\tau\in\overline\Delta$ and
$\eta\in\pl\Delta,\ \eta\not=\tau$, we say that a generator
$f\in\gtp$ belongs to the subcone $\gtep$ if it vanishes at the
point $\eta$, i.e., $\angle\lim\limits_{z\to\eta}f(z)=0$ and the
angular derivative at the point $\eta$
\[
f'(\eta):=\angle\lim_{z\to\eta}\frac{f(z)}{z-\eta}
\]
exists finitely.

\noindent$\bullet$ We say that a function $h\in\Spt$ ($h\in\Stt$)
belongs to the subclass $\Spte$ ($\Stte$) if the angular limit
\[
Q_h(\eta):=\angle\lim_{z\to\eta}\frac{(z-\eta)h'(z)}{h(z)}
\]
exists finitely and is different from zero.

\vspace{2mm}

\noindent{\bf Remark 3.} We recall that if $\zeta\in\pl\Delta$ and
$g\in\Hol(\Delta,\C)$ is such that
$\angle\lim\limits_{z\to\zeta}g(z)=:g(\zeta)$ exists finitely, the
expression
\[
Q_g(\zeta,z):=\frac{(z-\zeta)g'(z)}{g(z)-g(\zeta)}
\]
is called the {\it Visser--Ostrowski quotient} of $g$ at $\zeta$
(see \cite{PC1}). If for some ${h\in\Hol(\Delta,\C)}$ we have
$\angle\lim\limits_{z\to\zeta}h(z)=\infty$, then the
Visser--Ostrowski quotient of $h$ is defined by
\[
Q_h(\zeta,z):=Q_{1/h}(\zeta,z).
\]
A function $g$ is said to satisfy the {\it Visser--Ostrowski
condition} if
\[
Q_g(\zeta):=\angle\lim_{z\to\zeta}Q_g(\zeta,z)=1.
\]

In this context, we recall also that $g\in\Hol(\Delta,C)$ is
called {\it conformal at} $\zeta\in\pl\Delta$ if the angular
derivative $g'(\zeta)$ exists and is neither zero nor infinity;
$g$ is called {\it isogonal at $\zeta$} if the limit of
$\dst\arg\,\frac{g(z)-g(\zeta)}{z-\zeta}\,$ as $z\to\zeta$ exists.

It is clear that any function $g$ conformal at a boundary point
$\zeta$ is isogonal at this point. Also, it is known (see
\cite{PC1}) that any function $g$ isogonal at a boundary point
$\zeta$ satisfies the Visser--Ostrowski condition at this point,
i.e., $Q_g(\zeta)=1$.

So it is natural to say that $g$ satisfies a {\it generalized
Visser--Ostrowski condition} if
$Q_g(\zeta):=\angle\lim\limits_{z\to\zeta}Q_g(\zeta,z)$ exists
finitely and is different from zero. Thus each function
$h\in\Spte$ ($h\in\Stte$) satisfies a generalized
Visser--Ostrowski condition at the boundary point $\eta$.

\vspace{3mm}

To proceed, we note that the inequality $\eta\not=\tau$ implies
that for each ${h\in\Spte}$
\[
\angle\lim_{z\to\eta}h(z)=\infty.
\]

The following fact is an immediate consequence of Lemma~\ref{1}.

\begin{lemma}\label{2}
Let $h\in\Spt$ and $f\in\gtp$ be connected by {\rm(\ref{stareq})}.
Then $h$ belongs to $\Spte$ if and only if $f\in\gtep$. In this
case,
\[
Q_h(\eta)=\frac\mu{f'(\eta)}\,.
\]
\end{lemma}

We require two representation formulas for the classes of starlike
functions $\Stt$ and $\Stte$. For a boundary point $w$, denote by
$\delta_w$ the Dirac measure at this point.

\begin{lemma}[cf. \cite{GAW} and \cite{E-S3}]\label{3}
Let $\tau\in\overline\Delta$ and $\eta\in\pl\Delta,\
\eta\not=\tau$. Let $h\in\Hol(\Delta,\C)$ satisfy $h(\tau)=0$.
Then

(i) $h\in\Stt$ if and only if it has the form
\be\label{int0}
h(z)=C(z-\tau)(1-z\bar\tau)\cdot
\exp\left[-2\oint\limits_{\pl\Delta}
\log(1-z\bar\zeta)d\widetilde{\sigma}(\zeta) \right],
\ee
where $d\widetilde{\sigma}$ is an arbitrary probability measure on
the unit circle and $C\not=0$.

(ii) Moreover, $h\in\Stte$ if and only if it has the form
\be\label{int}
h(z)=&&C(z-\tau)(1-z\bar\tau)(1-z\bar\eta)^{-2a}\cdot \nonumber\\
&&\cdot\exp\left[-2(1-a)\oint\limits_{\pl\Delta}
\log(1-z\bar\zeta)d\sigma(\zeta) \right],
\ee
where $d\sigma$ is a probability measure on the unit circle
singular relative to $\delta_\eta$, $C\not=0$ and $a\in(0,1]$. In
this case, $Q_h(\eta)=-2a$.
\end{lemma}

\noindent{\bf Remark 4.} The constant $C$ can be chosen starting
from a normalization of functions under consideration. On the
other hand, since a starlike function $h$ is a solution of a
linear homogeneous equation (see (\ref{stareq})), $C$ arises in
the integration process of this equation.

\pr First, suppose that $\tau=0$, and let $h\in\Hol(\Delta,\C)$ be
normalized by ${h(0)=0}$ and $h'(0)=1$. A well-known criterion of
R. Nevanlinna asserts that $h\in\Star[0]$ if and only if
\[
q(z):=\frac{zh'(z)}{h(z)}
\]
has positive real part. (Note that the same fact follows by
(\ref{stareq}), because by the Berkson--Porta representation
formula (\ref{3a}), a generator $f\in\G[0]$ has the form
$f(z)=zp(z)$ with $\Re p(z)>0$).)

Representing $q$ by the Riesz--Herglotz formula, we write
\[
\frac{zh'(z)}{h(z)} =\oint\limits_{\pl\Delta}
\frac{1+z\bar\zeta}{1-z\bar\zeta}\, d\widetilde{\sigma}(\zeta)
\]
with some probability measure $d\widetilde{\sigma}$. Integrating
this equality, we get
\be\label{inta}
h(z)=z\exp\left[-2\oint_{\pl\Delta}\log(1-z\bar\zeta)
d\widetilde{\sigma}(\zeta)\right].
\ee
So we have proved (\ref{int0}) for the case $\tau=0$.

Now let $\tau\in\Delta$ be different from zero, and suppose
$h(\tau)=0$. It was proved by Hummel (see \cite{HJA-68},
\cite{HJA-69} and \cite{SD}) that $h\in\Stt$ if and only if
\[
\frac{z}{(z-\tau)(1-z\bar\tau)}\,h(z)\in\Star[0].
\]
Thus, (\ref{inta}) implies (\ref{int0}) for the interior location
of $\tau$. The reverse consideration and Hummel's criterion show
that if $h$ satisfies (\ref{int0}) with $\tau\in\Delta$, it must
be starlike.

Finally, let $\tau\in\pl\Delta$. Following Lyzzaik \cite{LA} (see
also \cite{E-G-R-S}) one can approximate $h\in\Stt$ by a sequence
$\{h_n\}$ of functions starlike with respect to those interior
points $\tau_n$ which converge to $\tau$. Also, one can assume
that $h_n(0)=h(0)$. Representing each function $h_n$ by
(\ref{int0})
\[
h_n(z)=C_n(z-\tau_n)(1-z\bar\tau_n)\cdot
\exp\left[-2\oint\limits_{\pl\Delta}
\log(1-z\bar\zeta)d\widetilde{\sigma_n}(\zeta) \right],
\]
we see that
\[
h(0)=h_n(0)=-C_n\tau_n.
\]
Thus $C_n\to-\frac{h(0)}\tau$. Since the set of all probability
measures is compact, $\{d\widetilde{\sigma_n}\}$ has a subsequence
converging to some probability measure $d\widetilde{\sigma}$.
Therefore, any function $h\in\Stt$ has the form (\ref{int0}).

To prove the converse assertion, we suppose that $h$ has the form
(\ref{int0}) with $\tau\in\pl\Delta$. Note that $h$ is starlike if
and only if the function $ah(cz),\ {a\neq0},\ {|c|=1}$, is.
Therefore, without loss of generality, one can assume that $h$ is
normalized by $h(0)=1$, i.e.,
\[
h(z)=(1-z)^2\cdot \exp\left[-2\oint\limits_{\pl\Delta}
\log(1-z\bar\zeta)d\widetilde{\sigma}(\zeta) \right].
\]
Differentiating the latter formula, one sees that $h$ satisfies a
modified Robertson inequality (see \cite{SH-SEM} and
\cite{E-R-S1})
\be\label{gri}
\Re\left[ \frac{zh'(z)}{h(z)}+\frac{1+z}{1-z} \right]>0.
\ee
A main result of \cite{SH-SEM} and Theorem 7 \cite{E-R-S1} imply
that $h$ is a starlike function with respect to a boundary point
with $h(1)=0$, i.e., $h\in\Star[1]$. The first assertion is
proved.

Let
\[
d\widetilde{\sigma}=a\delta_\eta+(1-a)d\sigma,\quad 0\le b\le1,
\]
be the Lebesgue decomposition of $d\widetilde{\sigma}$ relative to
the Dirac measure $\delta_\eta$, where the probability measures
$d\sigma$ and $\delta_\eta$ are mutually singular. Using this
decomposition, we rewrite (\ref{int0}) in the form (\ref{int}).

Now we calculate
\be\label{p}
&&Q_h(\eta)=\angle\lim_{z\to\eta}
\frac{h'(z)(z-\eta)}{h(z)}\nonumber \\ && =
\angle\lim_{z\to\eta}(z-\eta) \left[
\frac{\left((z-\tau)(1-z\bar\tau)\right)'}{(z-\tau)(1-z\bar\tau)}
+ \frac{2a\bar\eta}{1-z\bar\eta}
+2(1-a)\oint_{\pl\Delta}\frac{\bar\zeta}{1-z\bar\zeta}
\,d\sigma(\zeta) \right]\nonumber
\\ && = -2a+ 2(1-a)\, \angle\lim_{z\to-1}
\oint_{\pl\Delta}\frac{\bar\zeta(z-\eta)}{1-z\bar\zeta}\,d\sigma(\zeta)
\ee

Noting that
\[
\left| \frac{\bar\zeta(z-\eta)}{1-z\bar\zeta} \right|
\le\frac{|z-\eta|}{1-|z|}\,,
\]
we see that the integrand in the last expression of (\ref{p}) is
bounded on each nontangential approach region
$D_{k,\eta}:=\left\{z:\,|z-\eta|<k(1-|z|)\right\},\ k\ge1,$ at the
point $\eta$. Since the measures $d\sigma$ and $\delta_\eta$ are
mutually singular, we conclude by the Lebesgue convergence theorem
that the last integral in (\ref{p}) is equal to zero, so
\[
Q_h(\eta)=-2a.
\]
The proof is complete. \epr

The following results are angle distortion theorems for starlike
and spirallike functions of the classes $\Stte$ and $\Spte$
respectively.

\begin{lemma}[cf. \cite{SS} and \cite{E-S3}]\label{4}
Let $h\in\Stte$ with $Q_h(\eta)=\nu$. Denote
\be\label{aux2}
\theta:=\lim_{r\to1^-}\arg h(r\eta).
\ee
Then the image $h(\Delta)$ contains the wedge
\be\label{wedge2}
W=\left\{w\in\C:\ |\arg w -\theta|<\frac{|\nu|\pi}2\right\}
\ee
and contains no larger wedge with the same bisector.
\end{lemma}

\pr By Lemma~\ref{3}, the function $h$ has the form (\ref{int})
with $\nu=-2a$.

First we show that the image $h(\Delta)$ contains the wedge $W$
defined by~(\ref{wedge2}).

Since (as mentioned above)
$\angle\lim\limits_{z\to\eta}h(z)=\infty$, for each
$\delta\in\left(0,\frac\pi2\right)$ and each $R>0$, there exists
$r>0$ such that
\be\label{aux1}
|h(z)|>R
\ee
whenever
\[
z\in D_{r,\delta}:=\left\{z\in\Delta:\ |1-z\bar\eta|\le r,\
|\arg(1-z\bar\eta)|\le\delta \right\}.
\]
Lemma~\ref{3} and the Lebesgue convergence theorem imply the
existence of
\bep
&&\hspace{-2cm}\lim_{z\to\eta}\arg\frac{h(z)}{(1-z\bar\eta)^\nu}\hfill\\
&&\hspace{-2cm}=\arg\Bigl(C(\eta-\tau)(1-\eta\bar\tau)\Bigr)
-2(1-a)\lim_{z\to\eta}
\oint\limits_{\pl\Delta}\arg(1-z\bar\zeta)d\sigma(\zeta).
\eep
On the other hand, by formula (\ref{int}), we have
\bep
&&\hspace{-2cm}\theta=\lim_{r\to1^-}\arg h(r\eta)\\
&&\hspace{-2cm} =\arg\Bigl(C(\eta-\tau)(1-\eta\bar\tau)\Bigr) -
2(1-a) \lim_{r\to1^-}
\oint_{\pl\Delta}\arg(1-r\eta\bar\zeta)d\sigma(\zeta).
\eep
Therefore,
\[
\lim_{z\to\eta}\arg\frac{h(z)}{(1-z\bar\eta)^\nu}=\theta.
\]

Thus, decreasing $r$ (if necessary), we have
\[
\theta-\varepsilon< \arg\frac{h(z)}{(1-z\bar\eta)^\nu}<
\theta+\varepsilon
\]
for all $z\in D_{r,\delta}$. So, for each point $z$ belonging to
the arc
\[
\Gamma:=\left\{z\in\Delta:\ |1-z\bar\eta|= r,\
|\arg(1-z\bar\eta)|\le\delta \right\}\subset  D_{r,\delta},
\]
i.e., $z=\eta(1-re^{it}),\,|t|\le\delta$, we get
\[
\theta-\varepsilon -t|\nu|< \arg h(z) <\theta+\varepsilon -t|\nu|.
\]
In particular,
\be\label{aux3}
\arg h(\eta(1-re^{i\delta}))< \theta+\varepsilon -\delta|\nu|
\ee
and
\be\label{aux4}
\arg h(\eta(1-re^{-i\delta}))> \theta-\varepsilon +\delta|\nu|.
\ee

Thus, the curve $h(\Gamma)$ lies outside the disk $|z|\le R$ and
joins two points having arguments less than
$\theta+\varepsilon-\delta|\nu|$ and greater than
$\theta-\varepsilon +\delta|\nu|$, respectively. Since $h$ is
starlike, we see that $h(\Delta)$ contains the sector
\[
\left\{w\in\C:\ |w|<R,\ |\arg w-\theta|<\delta|\nu|-\varepsilon
\right\}.
\]
Since $R$ and $\varepsilon$ are arbitrary, one concludes
\[
\left\{w\in\C:\ |\arg w-\theta|<\delta|\nu|\right\}\subset
h(\Delta).
\]
Letting $\delta$ tend to $\frac\pi2$, we obtain
\[
W=\left\{w\in\C:\ |\arg w-\theta|<\frac{|\nu|\pi}2\right\}\subset
h(\Delta).
\]

Further, since $h$ is a starlike function, $\arg h(e^{i\varphi})$
is an increasing function in
$\varphi\in(\arg\eta-\pi,\arg\eta+\pi)$. So the limits
\[
\lim_{\varphi\to(\arg\eta)^\pm}\arg h(e^{i\varphi})
\]
exist. Let $\varphi_{n,+}\to(\arg\eta)^+$ and
$\varphi_{n,-}\to(\arg\eta)^-$ be two sequences such that the
values $h(e^{i\varphi_{n,\pm}})$ are finite. Then, once again by
Lemma~\ref{3},
\bep
&&\lim_{n\to\infty}\arg h(e^{i\varphi_{n,+}})-\arg
h(e^{i\varphi_{n,-}})    \\&&= \lim_{n\to\infty}
\left((\arg(1-e^{i\varphi_{n,+}}\bar\eta))^\nu -
(\arg(1-e^{i\varphi_{n,-}}\bar\eta))^\nu \right) =|\nu|\pi.
\eep
Therefore, the image contains no wedge of angle larger than
$|\nu|\pi$. Thus, the wedge $W$ defined by (\ref{wedge2}) is the
largest one contained in $ h(\Delta)$.

The proof is complete. \epr

Let $\lambda\in\Lambda=\{w \in\C:\ |w-1|\leq 1,\ w\neq 0\}$ and
$\theta\in[0,2\pi)$ be given. Define the function
$h_{\lambda,\theta}\in\Hol(\Delta)$ by
\begin{equation}\label{1c}
h_{\lambda,\theta }(z)=e^{i\theta}\left( \frac{1-z}{1+z}\right)
^\lambda.
\end{equation}
Here and in the sequel, we choose a single-valued branch of the
analytic function $w^\lambda$ such that $1^\lambda=1.$

\begin{defin}\label{def4}
The set $W_{\lambda,\theta }=h_{\lambda,\theta }(\Delta)$ is
called a canonical $\lambda$-spiral wedge with midline
${l_{\theta,\lambda}=\{w\in\C:\ w=e^{i\theta+t\lambda},\
t\in\R\}.}$
\end{defin}

To explain this definition, let us observe that
$h=h_{\lambda,\theta}$ is a solution of the differential equation
\[
\lambda h(z)=h'(z) f(z)
\]
normalized by the conditions $h(0)=e^{i\theta},\ h(1)=0,$ where
$f$ is given by
\[
f(z)=\frac{1}{2}\,(z^{2}-1).
\]
Since $f\in\gup$ with $f'(1)=1$ and $\lambda\in\Lambda$, it
follows by Lemma~1 that $h$ is a $\lambda$-spirallike function
with respect to the boundary point $h(1)=0$. Moreover, $f$ is a
generator of a one-parameter group (flow) of hyperbolic
automorphisms of $\Delta$ having two boundary fixed points $z=1$
and $z=-1.$ Hence, for each $w\in W_{\lambda,\theta}$ and
$t\in\R=(-\infty,\infty)$, the spiral curve $e^{-t\lambda}w$
belongs to $W_{\lambda,\theta }$ (see~(\ref{schroder})).

In \cite{A-E-S}, the notion of ``angle measure" for spirallike
domains with respect to a boundary point was introduced. It can be
shown that a $\lambda$-spiral wedge is of angle measure
$\pi\lambda$.

Finally, we see that for real $\lambda\in(0,2]$, the set
$W_{\lambda,\theta}$ is a straight wedge (sector) of angle
$\pi\lambda$, whose bisector is $l_\theta=\{w\in\C:\ \arg
w=\theta\}.$

\begin{lemma}\label{5}
Let $h\in\Spt$ be a $\mu$-spirallike function on $\Delta$. Then
the image $h(\Delta)$ contains a canonical $\lambda$-spiral wedge
with
\be\label{l51}
\arg\lambda=\arg\mu
\ee
if and only if $h\in\Spte$ for some $\eta\in\pl\Delta$. Moreover,
if $Q_h(\eta)=\nu$, then the canonical wedge
$W_{-\nu,\theta}\subset h(\Delta)$ for some $\theta\in[0,2\pi)$;
and it is maximal in the sense that there is no spiral wedge
$W_{\lambda,\theta}\subset h(\Delta)$ with $\lambda$ satisfying
{\rm(\ref{l51})} which contains $W_{-\nu,\theta}$ properly.
\end{lemma}

\pr First, given $h\in\Spte$ we construct $h_1\in\Spiral[1,-1]$
which is spirallike with respect to a boundary point whose image
eventually coincides with $h(\Delta)$ at $\infty.$ If
$\tau\in\pl\Delta$, we just set $h_{1}=h(\Phi(z)),$ where
$\Phi\in\Aut(\Delta)$ is an automorphism of $\Delta$ such that
$\Phi(1)=\tau$ and $\Phi(-1)=\eta.$

If $\tau\in\Delta$, we take any two points $z_{1}=e^{i\theta_{1}}$
and $z_{2}=e^{i\theta_{2}}$ such that ${w_{1}=h(z_{1})}$ and
$w_{2}=h(z_{2})$ exist finitely and
$\theta_1\in(\arg\eta-\epsilon,\arg\eta)$,\linebreak
${\theta_2\in(\arg\eta ,\arg\eta-\epsilon)}$, so the arc
$(\theta_{1,}\theta_{2})$ on the unit circle contains the
point~$\eta$.

Since $h$ is spirallike with respect to an interior point, it
satisfies the equation
\begin{equation}
\beta h(z)=h^{\prime }(z)f(z),  \label{3c}
\end{equation}
where $f\in\gtp$ and $\beta =f'(\tau )$, so $\arg\mu=\arg\beta$.
This means that for each $w\in h(\Delta)$ the spiral curve
$\{e^{-t\beta }w,\ t\geq 0\}$ belongs to $h(\Delta ).$ In turn,
the curves $l_{1}=\{z=h^{-1}(e^{-t\beta }w_{1}),\ t\geq 0\}$ and
$l_{2}=\{z=h^{-1}(e^{-t\beta }w_{2}),\ t\geq 0\}$ lie in
$\overline{\Delta }$ with ends in $z_{1}$\ and $\tau $ and $z_{2}$
and $\tau $, respectively.

Since $z_{1}\neq z_{2}$ and the interior points of $l_{1}$ and
$l_{2}$ are semigroup trajectories in $\Delta$, these curves do
not intersect except at their common end point $z=\tau.$
Consequently, the domain $D$ bounded by $l_{1},l_{2}$ and the arc
$(\theta _{1,}\theta _{2})$ is simply connected, and there is a
conformal mapping $\Phi$ of $\Delta$ such that $\Phi(\Delta)=D$
and $\Phi(-1)=\eta,\ \Phi(1)=\tau$. Now define $h_{1}(z)=h(\Phi
(z))$. It follows by our construction that $h_{1}(\Delta )\subset
h(\Delta)$ and $h_1$ is spirallike with respect to a boundary
point $h_1(1)=0.$ In addition, since $\Phi$ is conformal at the
point $z=-1$, it satisfies the Visser--Ostrowskii condition and we
have
\begin{eqnarray}\label{4c}
&&\angle\lim_{z\to-1}\frac{(z+1)h_1'(z)}{h_1(z)}
=\angle\lim_{z\to-1}\frac{(z+1)h'(\Phi (z))(\Phi
(z)-\eta)}{h(\Phi(z))(\Phi(z)-\eta)} \nonumber  \\
&=&\angle\lim_{z\to-1}\frac{(z+1)\Phi'(z)}{\Phi(z)-\Phi(-1)}\cdot
\angle \lim_{z\to-1}\frac{(\Phi (z)-\eta)h'(\Phi(z))}{h(\Phi (z))}
\nonumber \\
&=&\angle\lim_{z\to-1}\frac{(\Phi(z)+1)h'(\Phi(z))}{h(\Phi
(z))}\,.
\end{eqnarray}

Note also that $\Phi$ is a self-mapping of $\Delta $ mapping the
point $z=-1$ to $\eta$ and having a finite derivative at this
point.

It follows by the Julia--Carath\'eodory theorem, (see, for
example, \cite{SD}) that if $z$ converges to $-1$ nontangentially,
then $\Phi(z)$ converges nontangentially to $\eta=\Phi(-1)$. Then
(\ref{4c}) implies that
\begin{equation}
Q_{h_{1}}(-1)=\angle \lim_{z\to-1}\frac{(z+1)h_{1}'(z)}{h_{1}(z)}
\label{5c}
\end{equation}
exists finitely if and only if $h\in\Spte$ and
\begin{equation} \label{6c}
Q_{h_{1}}(-1)=Q_{h}(\eta).
\end{equation}

We claim that this last relation implies that $h_{1}(\Delta )$
contains a $(-\nu)$-spiral wedge $W_{-\nu,\theta}$ for some
$\theta\in[0,2\pi).$

To this end, observe that $h_{1}$ satisfies the equation
\[
\beta h_{1}(z)=h_{1}^{\prime }(z)\cdot f_{1}(z),
\]
where $f_1(z)=\frac{f(\Phi(z))}{\Phi'(z)}$ is a generator of a
semigroup of $\Delta $ with $f_1(1)=0$ and $f_1'(1)=\beta_1$ for
some $\beta_1>0$ such that
\[
\left\vert \beta -\beta_{1}\right\vert \leq \beta _{1}.
\]

Therefore, $h_{1}$ is a complex power of the function
$h_{2}\in\Hol(\Delta,\C)$ defined by the equation
\begin{equation}\label{7c}
\beta_1 h_{2}(z)=h_{2}'(z)f_{1}(z),\quad h_{2}(1)=0,
\end{equation}
i.e.,
\begin{equation}
h_{1}(z)=h_{2}^{\mu }(z),  \label{8c}
\end{equation}
where $\mu =\frac{\beta}{\beta_1}\neq 0,\ |\mu -1|\leq 1$, hence
$\arg\mu=\arg\beta$.

On the other hand, if we normalize $h_{1}$ by
$h_1^{1/\mu}(0)=h_2(0)$, equation~(\ref{8c}) has a unique solution
which is a starlike function with respect to a boundary point
$(h_{2}(1)=0).$ Obviously,
\be\label{l52}
Q_{h_2}(-1)=\frac{1}{\mu}\,Q_{h_1}(-1)\left(=\frac{1}{\mu}\,
Q_h(\eta)\right).
\ee

Note that $\nu _2:=Q_{h_2}(-1)$ is a negative real number, while
${\nu _1:=Q_{h_1}(-1)=\nu _{2}\mu}$ is complex.

Now it follows by Lemma \ref{4} that the starlike set
$h_2(\Delta)$ contains a straight wedge (sector) of a nonzero
angle $\sigma\pi$ for each $\sigma\in(0,|\nu_2|\pi]$. So the
maximal (straight) wedge $W\subset h_{2}(\Delta )$ is of the form
\[
W=W_{-\nu_2,\theta_2}=\left\{w\in\C:\
w=e^{i\theta_2}\left(\frac{1-z}{1+z}\right)^{-\nu_2}\right\},
\]
with
\begin{eqnarray*}
\theta_2 &=&\lim_{r\to 1^-}\arg h_2(-r)=\lim_{r\to 1^-}\arg
h_1^{\nu_2/\nu_1}(-r)= \\ &=& \nu_2\cdot\lim_{r\to 1^-}\arg
h_1^{1/\nu_1}(-r)= \nu_2\theta_1,
\end{eqnarray*}
where
\[
\theta_{1}=\lim_{r\rightarrow 1^{-}}\arg h_{1}^{1/\nu _{1}}(-r).
\]

Writing $W$ in the form
\[
W=\left\{ e^{i\varsigma}e^{t},\ t\in\R,\
\varsigma\in\left(\theta_2+\frac{\pi\nu_2}{2},\theta_2-\frac{\pi\nu_2}{2}\right)
\right\}
\]
and setting $\varsigma_1=\varsigma/\nu_2,\ s=t/\nu_2$, we see that
the set
\[
K:=W^\mu= \left\{e^{i\varsigma_1\nu_1}e^{s\nu_1},\ s\in\R,\
\varsigma_1\in\left(\frac{\theta_2}{\nu_2}
-\frac\pi2,\frac{\theta_2}{\nu_2}+\frac\pi2\right) \right\}
\]
is contained in $h_{1}(\Delta )$; hence in $h(\Delta).$ But
$\theta_2/\nu_2=\theta_1$ and $\nu_1=\nu\left(=Q_h(-1)\right)$;
hence $K$ is of the form
\bep
K=\left\{e^{i\varsigma_1\nu}e^{s\nu},\ s\in\R,\
\varsigma_1\in\left(\theta_1-\frac\pi2,\theta_1+\frac\pi2\right)\right\}\\
=\left\{e^{i\theta_1\nu}e^{i\varsigma_1\nu}e^{s\nu},\ s\in\R,\
\varsigma_1\in\left(-\frac\pi2,+\frac\pi2\right)\right\}.
\eep

Setting $\theta=\frac{|\nu|^2\theta_1}{\Re\nu}\in\R$, we get
\bep
i\theta_1\nu+s\nu=i\theta
+\nu\left(\frac{\theta\Re\nu}{|\nu|^2}+s-\frac{i\theta}\nu \right)
=i\theta +\nu\left(s-\frac{\theta\Im\nu}{|\nu|^2} \right).
\eep
Since $s$ takes all real values, so does
$t=s-\frac{\theta\Im\nu}{|\nu|^2}$. Therefore, the set $K$ has the
form
\bep
K=e^{i\theta}\left\{e^{i\varsigma_1\nu}e^{t\nu},\ t\in\R,\
\varsigma_1\in\left(-\frac\pi2\,,\,\frac\pi2\right)\right\},
\eep
i.e., coincides with $W_{-\nu,\theta}$. Finally, it follows by
(\ref{l52}) that $\lambda:=-\nu=|\nu_2|\mu$. This implies
(\ref{l51}).

Conversely, let $h$ be a $\mu$-spirallike function on $\Delta$
such that $h(\Delta)$ contains a canonical $\lambda$-spiral wedge
$W_{\lambda,\theta}$ for some $\lambda$ satisfying (\ref{l51}) and
$\theta\in[0,2\pi)$. Then for each $w_0\in W_{\lambda,\theta}$,
the curve $l:=\left\{w\in\C:\,w=e^{-t\lambda}w_0,\,t\in\R
\right\}$ belongs to $h(\Delta)$. Hence the curve
$h^{-1}(l)\subset\Delta$ joints the point $\tau\in\overline\Delta$
with a point $\eta\in\pl\Delta$. Again, as in the first step of
the proof, one can find a conformal mapping $\Phi\in\Hol(\Delta)$
with $\Phi(1)=\tau,\ \Phi(-1)=\eta$ such that $h_1=h\circ\Phi$ is
a $\mu$-spirallike function with respect to a boundary point
$h_1(1)=0$ and
\be\label{l53}
W_{\lambda,\theta}\subset h_1(\Delta)\subset h(\Delta).
\ee
Again the function $h_2=h_1^{1/\mu}$ is starlike with respect to a
boundary point, and $h_2(\Delta)$ contains the set
\[
K=\left\{w\in\C:\ w=e^{i\frac\theta\mu}\left(\frac{1-z}{1+z}
\right)^{\frac\lambda\mu} \right\}
\]
because of (\ref{l53}).

Setting $\,\dst\frac\lambda\mu=\kappa$ and
$\dst\theta_1=\theta\,\frac{\Re\mu}{|\mu|^2}\,$, we see by
(\ref{l51}) that $\kappa$ is real and $K$ can be written as
\[
K=\left\{w\in\C:\ w=Re^{i\theta_1}\left(\frac{1-z}{1+z}
\right)^\kappa \right\},
\]
with $R=\exp\left[\frac{\theta_1\Im\mu}{\Re\mu}\right]$ real and
positive.

Hence, $h_2(\Delta)$ contains a straight canonical wedge\linebreak
$W_{\kappa,\theta_1}=\left\{w\in\C:\
w=e^{i\theta_1}\left(\frac{1-z}{1+z} \right)^\kappa \right\}$ with
$0<\kappa|\nu_2|$, where $\nu_2=Q_{h_2}(-1)$ exists finitely and
$W_{|\nu_2|,\theta_1}$ is the maximal wedge contained in
$h_2(\Delta)$. But, as before, we have
\[
\nu=Q_h(\eta)=\mu Q_{h_2}(-1)=\mu\nu_2.
\]
The latter relations show that $\nu$ is finite and $\lambda$ must
satisfy the conditions $\arg\lambda=\arg\mu=\arg(-\nu)$ and
$0<|\lambda|\le|\nu|$. So the wedge $W_{-\nu,\theta}$ is a maximal
wedge contained in $h(\Delta)$ satisfying condition (\ref{l51}).
The lemma is proved. \epr

\noindent{\bf Remark 5.} By using Lemma \ref{4} and the proof of
Lemma \ref{5}, one can show that the number $\theta$ in the
formulation of Lemma 5 is defined by the formula
\[
\theta=\frac{|\nu|^2}{\Re\nu}\,\lim_{r\to1^-}\arg h^{1/\nu}(-r).
\]
For real $r$, this formula coincides with (\ref{aux2}). Hence, in
fact, Lemma \ref{5} contains Lemma \ref{4}.

\vspace{3mm}

\begin{lemma}
Let $f\in\gtep$ for some $\tau\in\overline\Delta$ (which is the
Denjoy--Wolff point for the semiflow $\Ss$ generated by $f$) with
$\beta=f'(\tau)>0$ and some $\eta\in\pl\Delta$, such that
$f(\eta):=\angle\lim\limits_{z\to\eta}f(z)=0$ and
$\gamma=f'(\eta)=\angle\lim\limits_{z\to\eta}f'(z)$ exists
finitely.

The following assertions hold.

(i) If $\tau\in\Delta$, then $\gamma <-\frac{1}{2}\Re\beta$.

(ii) If $\tau\in\pl\Delta$, then $\gamma\leq-\beta<0$ and the
equality $\gamma=-\beta$ holds if and only if
$f\subset\aut(\Delta)$ or, what is the same,
$\Ss\subset\Aut(\Delta )$ consists of hyperbolic automorphisms of
$\Delta$.
\end{lemma}

\pr (i) Let $\tau\in\Delta$. Then $f\in\gtp$ admits the
representation
\[
f(z)=(z-\tau)(1-z\bar\tau)p(z)
\]
with $\Re p(z)>0,\ z\in \Delta$ and
\[
\beta\left(=f'(\tau)\right)=(1-|\tau|^2)p(\tau ).
\]
Assume that for some $\eta\in\pl\Delta$
\[
f(\eta ):=\angle \lim_{z\rightarrow \eta }f(z)=0
\]
and
\[
\gamma =\angle\lim_{z\rightarrow \eta }\frac{f(z)}{z-\eta }
\]
exists finitely. Then $\angle\lim_{z\rightarrow \eta }p(z)=0$, and
\[
\gamma=\eta |\eta-\tau|^2\cdot p'(\eta ),
\]
where
\[
p'(\eta )=\angle \lim_{z\rightarrow \eta }\frac{p(z)}{z-\eta }.
\]

To find an estimate for $p'(\eta)$, we introduce a function $p_1$
of positive real part by the formula
\[
p_1(z)=(1-\left\vert \tau \right\vert ^{2})p(m(z)),
\]
where
\[
m(z)=\frac{\tau -z}{1-z\bar\tau}
\]
is the M\"obius transformation (involution) taking $\tau$ to $0$
and $0$ to $\tau$. Thus
\[
p_1(0)=(1-\left\vert \tau \right\vert ^{2})p(\tau )=\beta;
\]
and, setting  $\eta _{1}=m(\eta )$, we have
\bep
p_1'(\eta_1)&=&(1-|\tau|^2)p'(\eta)\cdot m'(\eta_1) \\
&=&\frac{1-|\tau|^2}{m'(\eta)}\cdot
p'(\eta)=-(1-\eta\bar\tau)^2p'(\eta ).
\eep

On the other hand, using the Riesz--Herglotz formula for the
function $q=1/p\,$ we obtain
\bep
\angle\lim_{z\to\eta_1}(z-\eta_1)q(z)
&&=\angle\lim_{z\to\eta_1} \int_{\pl\Delta} \frac{(z-\eta_1)(1+z\bar\zeta)}{1-z\bar\zeta}d\mu_q(\zeta)= \\
&&=-\eta_1\cdot\angle\lim_{z\to\eta_1}\int_{\pl\Delta}
\frac{(1-z\bar\eta_1)(1+z\bar\zeta)}{1-z\bar\zeta}d\mu _q(\zeta)= \\
&&=-\eta_12\mu_q(\eta_1),
\eep
where $\mu_q$ is a positive measure on $\pl\Delta$ such that
$\int_{\pl\Delta}d\mu_q(\zeta)=\Re q(0)$.\linebreak Consequently,
\bep
p_1'(\eta_1)&=&\angle\lim_{z\to\eta_1}\frac{p_1(z)}{z-\eta_1}=
\angle\lim_{z\to\eta_1}\frac{1}{(z-\eta_1)q(z)}= \\
&=&\frac{-\overline{\eta_1}}{2\mu_q(\eta_1)}=-(1-\eta\bar\tau)^2p'(\eta).
\eep
Hence
\[
p'(\eta)=\frac{\overline{\eta_1}}{(1-\eta\bar\tau)^2
2\mu_q(\eta_1)}
\]
and
\[
\gamma=\frac12\frac{\eta|\eta-\tau|^2\overline{\eta_1}}{(1-\eta\bar\tau)^2}\cdot
\frac{1}{\mu _{q}(\eta _{1})}.
\]

Since $\dst\mu_q(\eta_1)\leq\Re q(0)\leq\frac{1}{\Re
p_1(0)}=\frac{1}{\Re\beta}\,$, we have
\[
|\gamma|\geq \frac{1}{2}\Re\beta .
\]

Note that equality is impossible since otherwise $q$ (and hence
$p_1$ and $p$) are constant. But
$\angle\lim\limits_{z\to\eta_1}p(z)=0$, which means that
$p(z)\equiv 0.$

This proves assertion (i).

(ii) Let now $\tau\in\pl\Delta$. In this case, we know already
that
\[
\beta =f'(\tau)=\angle\lim_{z\to\tau }f'(z)>0.
\]

Without loss of generality, let us assume that $\tau=1$ and
$\eta=-1$. In other words, we assume that $f\in\G^+[1,-1]$. We
have to show that in that case
$\gamma=\angle\lim\limits_{z\to-1}f'(z)\le-\beta$, and equality
holds if and only if $f$ is a complete vector field.

Indeed, suppose to the contrary that $\gamma\in(-\beta ,0).$ Then
the function $g\in\Hol(\Delta,\C)$ defined by
\[
g(z)=f(z)+\frac{\gamma }{2}(z^{2}-1)
\]
belongs to the class $\G^+[1,-1],$ because this class is a real
cone. In addition,
\[
g'(1)=\beta +\gamma \geq 0,
\]
while
\[
g'(-1)=\gamma -\gamma =0.
\]
Then either $g(z)\equiv 0$, or $g\neq 0$ and both points $1$ and
$-1$ are sink points of the semigroup generated by $f$, which is
impossible. This contradiction shows that $g$ must be identically
zero, hence $\gamma =-\beta $ and
\[
f(z)=-\frac{\gamma }{2}(z^{2}-1).
\]

Thus $f$ belongs to $\aut(\Delta)$,  and the flow
$\Ss=\{F_{t}\}_{t\in\R}$ consists of hyperbolic automorphisms of
$\Delta .$ The lemma is proved. \epr

Now we are ready to prove our theorems. Since Theorem~2 is a
compliment of Theorem~1, we give their proofs simultaneously.

\noindent\textbf{Proof of Theorems 1 and 2.} We prove implications
$(i)\Longrightarrow(ii)\Longrightarrow(iii)\Longrightarrow(i)$ of
Theorem~1 successively, while assertions (a), (b) and (c) of
Theorem~2 will be obtained in the process. Let
$\Ss=\{F_{t}\}_{t\geq 0}$ be a semiflow on $\Delta$ generated by
$f\in\gtp$ with $\beta=f'(\tau)$, $\Re\beta>0$. Let
$h\in\Hol(\Delta,\C)$ be the associated spirallike (starlike)
function on $\Delta$ defined by equation (\ref{stareq}) with $\mu
=\beta.$ Then by Lemma~\ref{1}, $h$ satisfies Schr\"oder's
equation (\ref{schroder})
\be\label{schroder1} h(F_{1}(z))=e^{-t\beta }h(z) \ee for all
$t\geq 0$ and $z\in\Delta.$

{\bf Step 1 ((i)$\Longrightarrow $(ii)).} If now $f\in\gtep$ for
some $\eta\in\pl\Delta$, that is
$f(\eta)\left(=\angle\lim_{z\to\eta}f(z)\right)=0$ and
$\gamma=f'(\eta)\left(=\angle\lim_{z\to\eta}f'(z)\right)$ exists
finitely, then by Lemma~\ref{2} the function $h$ belongs to the
class $\Spte$ with
\[
Q_h(\eta)=\angle\lim_{z\to\eta}\frac{(z-\eta
)h'(z)}{h(z)}=\frac{\beta }{\gamma }\,.
\]
Since $\gamma \neq 0$ (actually, $\gamma <0$), $Q_h(\eta)$ is
finite.

In turn, Lemma \ref{5} implies that there is a non-empty (spiral)
wedge ${W\subset h(\Delta)}$ with vertex at the origin such that
for each $w\in W$ the spiral curve $\{e^{-t\beta }w\}$ belongs to
$W$, for all $t\in\R.$

Define the simply connected domain $\Omega\subset\Delta$ by
\[
\Omega =h^{-1}(W).
\]
Then the family $\widetilde{F_t}:\Omega\mapsto\Omega$
\[
\widetilde{F_t}(z)=h^{-1}\left( e^{-t\beta }h(z)\right) ,\quad
z\in \Omega ,\ t\in\R,
\]
forms a flow (one-parameter group) of holomorphic self-mappings of
$\Omega$. Comparing the latter formula with (\ref{schroder1}), we
see that for $t\geq 0$, $\widetilde{F_t}(z)=F_{t}(z)$ whenever
$z\in\Omega$ and
$\left(\left.F_t\right|_\Omega\right)^{-1}=\widetilde{F_{-t}}$.
Thus $\Ss\subset\Aut(\Omega).$

{\bf Step 2 ((ii)$\Longrightarrow $(iii)).} Let again
$\Ss=\{F_{t}\}_{t\geq 0}$ be a semiflow generated by $f\in\gtp$ so
that
\begin{equation}\label{1b}
\lim_{t\to\infty}F_{t}(z)=\tau\in\pl\Delta\quad\mbox{and}\quad
\Re\beta>0,\quad\mbox{where }\beta=f'(\tau ),
\end{equation}
and let $\Omega\subset\Delta$ be a simply connected domain such
that $\Ss\subset\Aut(\Omega).$ Let ${\psi:\Delta\mapsto\Omega}$ be
any Riemann conformal mapping of $\Delta$ onto $\Omega$. Consider
the flow $\{G_t\}_{t\in\R}\subset\Aut(\Delta)$ defined by
\begin{equation}\label{2b}
G_t(z)=\psi^{-1}(F_t(\psi(z))),\ t\in\R.
\end{equation}
In this case, $\psi$ is a conjugation for $G_{t}$ and $F_{t}$ for
each $t\in\R$, i.e.,
\begin{equation}\label{3b}
\psi(G_{t}(z))=F_{t}(\psi (z)),\quad z\in \Delta,\ t\in\R .
\end{equation}
Denote by $g\in\aut(\Delta)$ the generator of
$\{G_{t}\}_{t\in\R}$:
\[
g(z)=\dst\lim_{t\to0}\frac{z-G_{t}(z)}{t}\,.
\]
Then by (\ref{3b}), $\psi$ satisfies the differential equation
\begin{equation}
\psi'(z)\cdot g(z)=f(\psi (z)).  \label{4b}
\end{equation}

First we show that the family
$\{G_{t}\}_{t\in\R}\subset\Aut(\Delta)$ consists of hyperbolic
automorphisms or, what is the same, that it does not contains
neither elliptic nor parabolic automorphisms.

Indeed, suppose $\{G_{t}\}_{t\in\R}$ contains an elliptic
automorphism. Then there is a point $a\in\Delta$ such that
$G_{t}(a)=a$ for all $t\in\R$; hence $g(a)=0$ and $\Re g'(a)=0$.
By (\ref{4b}), $f(\psi(a))=0$; and thus $\psi(a)=\tau$. On the
other hand, differentiating (\ref{4b}) with respect to $z$ and
setting $z=a$, we get $g'(a)=f'(\tau)$. Hence $\Re f'(\tau)=0$,
which contradicts (\ref{1b}).

Thus $\{G_{t}\}$ has no interior fixed point in $\Delta$; hence
there are boundary points $\zeta_1$ and $\zeta_{2}$ such that
\begin{equation}\label{5b}
\lim_{t\rightarrow \infty }G_{t}(z)=\zeta _{1}\in\pl\Delta,\quad
z\in \Delta,
\end{equation}
and
\begin{equation}
\lim_{t\rightarrow -\infty }G_{t}(z)=\zeta _{2}\in\pl\Delta,\quad
z\in \Delta.   \label{6b}
\end{equation}
To show that the family $\{G_{t}\}_{t\in\R}$ does not contain a
parabolic automorphism it is sufficient to prove that $\zeta_1\neq
\zeta_2$.

To this end, we again consider the associated spirallike
(starlike) function $h$ defined by equation (\ref{stareq}) with
$\mu=\beta$ and normalized by the conditions $h(\tau)=0,\ h'(\tau
)=1$ if $\tau\in\Delta$ or by $h(\tau)=0$ and $h(0)=1$ if $\tau
\in\pl\Delta$ (see Lemma~\ref{1}). Define $h_0\in\Hol(\Delta,\C)$
by
\begin{equation}
h_{0}(z)=h(\psi (z)).  \label{7b}
\end{equation}
Since $h$ satisfies Schr\"oder's equation (\ref{schroder1}), it
follows from (\ref{3b}) that for all $t\geq 0$,
\begin{eqnarray*}
h_{0}(G_{t}(z)) &=&h(\psi (G_{t}(z))=h(F_{t}(\psi (z)) \\
&=&e^{-t\beta }h(\psi (z))=e^{-t\beta }h_{0}(z).
\end{eqnarray*}
Since the mapping $G_t\in\Hol(\Delta)$ is an automorphism of
$\Delta$ for each $t\in\R^+$, we have, in fact,
\begin{equation}
h_{0}(G_{t}(z))=e^{-t\beta }h_{0}(z)  \label{8b}
\end{equation}
for all $t\in\R$.

From (\ref{8b}) we conclude that $h$ is a univalent spirallike
(starlike) function on $\Delta$. Moreover, (\ref{8b}) and
Corollary 2.17 of \cite{PC1} imply that
\[
\angle \lim_{z\in \zeta _{1}}h_{0}(z)=0,
\]
while
\[
\angle \lim_{z\in \zeta _{2}}h_{0}(z)=\infty.
\]
Thus $\zeta _{1}\neq \zeta _{2}$, and it follows that
$\{G_{t}\}_{t\in\R}$ consists of hyperbolic automorphisms.

Now observe that $W=h_{0}(\Delta )$ is a spirallike (starlike)
wedge with vertex at the origin belonging to $h(\Delta)$. Since
all the points of $\pl h(\Delta)$ are admissible,
$\psi=h^{-1}\circ h_0$ is a homeomorphism of $\overline{\Delta}$
onto $\overline{\Omega}$; hence $\pl\Omega$ is a Jordan curve. Now
(\ref{3b}) implies that
\begin{equation}\label{9b}
\lim_{t\to\infty }\psi(G_t(z))=\lim_{t\to\infty
}F_t(\psi(z))=\tau
\end{equation}
and
\begin{equation}\label{10b}
\lim_{t\to-\infty }\psi(G_{t}(z))=\lim_{t\to-\infty
}F_{t}(\psi(z))=\eta
\end{equation}
for some $\eta\in\overline{\Delta}.$ Applying again Corollary 2.17
in \cite{PC1}, we obtain
\begin{equation}\label{11b}
\psi(\zeta_1):=\lim_{z\to\zeta_1}\psi(z)=\tau
\end{equation}
and
\begin{equation}\label{12b}
\psi (\zeta _{2}):=\angle \lim_{z\rightarrow \zeta _{2}}\psi
(z)=\eta.
\end{equation}

Thus $\eta\neq\tau$ and, moreover, $\eta\in\pl\Delta$. Indeed, if
$\eta$ is an interior point of $\Delta$, then
\[
\eta=\psi(\zeta_2)=\psi(G_t(\zeta_2))=F_t(\psi(\zeta_2))=
F_t(\eta),\quad t\geq 0,
\]
i.e., it must be an interior fixed point for all $F_t\in\Ss,\
t\ge0,$ which is impossible.

So $\tau\in\pl\Omega$ by (\ref{11b}), and
$\eta\in\pl\Delta\cap\pl\Omega$ by (\ref{12b}).

To show that equation (\ref{f**}) has a locally univalent (even
univalent) solution $\varphi\in\Hol(\Delta)$ for some $\alpha >0$,
we use a M\"obius transformation ${m\in\Aut(\Delta)}$ such that
$m(1)=\zeta_1$ and $m(-1)=\zeta_2$. Then $\psi_1=\psi\circ m$ is a
conformal mapping of $\Delta$ onto $\Omega$ with normalization
\[
\psi_1(1)=\tau ,\quad \psi _{1}(-1)=\eta .
\]

For $s\in (-1,1)$, define another conformal mapping $\varphi_s$ of
$\Delta$ onto $\Omega$ by
\[
\varphi_s(z):= \psi_1\left(\frac{z-s}{1-zs}\right),\quad -1<s<1.
\]
Clearly $\varphi _{s}(1)=\psi _{1}(1)=\tau $ and
$\varphi_{s}(-1)=\psi _{1}(-1)=\eta .$ Note also that
$l=\{z=\varphi_s(0)\left(=\psi_1(-s)),\ s\in[-1,1]\right)\}$ is a
continuous curve joining the points $z=1$ and $z=-1$, and so
$l_{1}=\{z=h(\varphi _{s}(0))(=h(\psi _1(-s)))\}$ is a continuous
curve joining $h(\tau )=0$ and $h(\eta )=\infty$. Hence, there
exists $s\in (-1,1)$ such that $|h(\varphi_s(0))|=1.$

Thus there exists a homeomorphism $\varphi(=\varphi_s)$ of
$\overline{\Delta}$ onto $\overline{\Omega}$ holomorphic in
$\Delta$ such that $\varphi(1)=\tau,\ \varphi(-1)=\eta$ and
$h(\varphi(0))=e^{i\theta}$ for some $\theta\in\R$.

Since the mapping $\psi$ in our previous consideration was
arbitrary, we can replace it by $\varphi.$ In this case, the
``new" flow $\{G_{t}\}_{t\in 0}$ defined by
\[
G_{t}(z)=\varphi ^{-1}(F_{t}(\varphi(z)))
\]
is a one-parameter group of hyperbolic automorphisms of $\Delta$
having the fixed points $z=1$ and $z=-1$ on $\pl\Delta$. In turn,
its generator $g\in\Hol(\Delta,\C)$ must have the form
\begin{equation}
g(z)=\frac{\alpha }{2}(z^{2}-1)\, ,  \label{13b}
\end{equation}
where $\alpha =g^{\prime }(1)>0.$

Hence, equation (\ref{4b}) (with $\varphi$ in place of $\psi$)
becomes (\ref{f**})
\begin{equation}\label{14b}
\alpha\varphi'(z)(z^{2}-1)=2f(\varphi (z)).
\end{equation}

Combining this with (\ref{stareq}), we show that $\alpha$ must be
greater than or equal to $-\gamma>0$. Namely, defining
$h_{0}\in\Spiral[1]$ as in (\ref{7b}) by
\begin{equation}\label{15-b}
h_0(z)=h(\varphi(z)),
\end{equation}
we have from (\ref{14b}) and (\ref{stareq}) that
\begin{equation}
\beta h_{0}(z)=\frac{\alpha }{2}(z^{2}-1)h_{0}^{\prime }(z)  \label{16-b}
\end{equation}
with $h_{0}(0)=h(\varphi (0))=e^{i\theta}$ for some
$\theta\in[0,2\pi)$. Solving this equation, we obtain
\be\label{aux5}
h_0(z)=e^{i\theta}\left(\frac{1-z}{1+z}\right)
^{\beta /\alpha }
\ee
with $\alpha =g'(1).$

On the other hand, by Lemma \ref{5}, the maximal (spiral) wedge
contained in $h(\Delta)$ is of the form
$W_{-\nu,\theta}=\left\{w\in\C:\
w=e^{i\theta}\left(\frac{1-z}{1+z}\right)^{-\nu} \right\}$, where
\be\label{aux6}
\nu =\angle\lim_{z\to\eta}\frac{(z-\eta
)h'(z)}{h(z)}=\angle\lim_{z\to\eta}\frac{z-\eta}{f(z)}=
\frac{\beta}{\gamma}\,.
\ee
Thus $\gamma$ is finite and $\varphi=h^{-1}\circ h_0$ is a
well-defined self-mapping of $\Delta$ if and only if
$\alpha\geq-\gamma$. This completes the proof  of the implication
((ii)$\Longrightarrow $(iii)) of Theorem~1, as well as assertions
(a) and (b) of Theorem~2. Note in passing that we have also proved
the implication (ii)$\Longrightarrow $(i) of Theorem~1.

{\bf Step 3 ((iii)$\Longrightarrow $(i)).} Suppose now that
$\varphi\in\Hol(\Delta)$ is locally univalent and satisfies
(\ref{f**}) for some $\alpha\in\R^+$. Solving this differential
equation explicitly, we get
\be\label{n1}
\alpha\int_{\varphi(0)}^{\varphi(z)}\frac{dw}{f(w)}=\int_0^z\frac{2dz}{z^2-1}=
\log\left(\frac{1-z}{1+z}\right).
\ee
Since $\varphi'(z)\neq0,\ z\in\Delta$, we have by (\ref{f**}) that
there is no $z\in\Delta$ such that $\varphi(z)=\tau$. So if $l$ is
a curve joining $0$ and $z$, the curve $\varphi(l)$ joining
$\varphi(0)$ and $\varphi(z)$ does not contain $\tau$.

Consider now the differential equation (\ref{stareq}) with initial
data ${h(\varphi(0))=1}$. Separating variables in this equation,
we see that
\be\label{n2}
\beta\int_{\varphi(0)}^{\varphi(z)}
\frac{dw}{f(w)}=\int_1^{h(\varphi(z))}\frac{dh}{h}=\log
(h(\varphi(z)).
\ee
Comparing (\ref{n1}) with (\ref{n2}), we have
\[
\log\left(h(\varphi(z)\right)=
\frac\beta\alpha\,\log\left(\frac{1-z}{1+z}\right),
\]
or
\[
h(\varphi(z))=\left(\frac{1-z}{1+z}\right)^{\frac\beta\alpha}.
\]
This equality implies that the set
$\dst\left\{\left(\frac{1-z}{1+z}\right)^{\frac\beta\alpha}:\
z\in\Delta\right\}$ is a subset of $h(\Delta)$, so this set is
different from $\C\setminus\{0\}$. It follows by \cite{A-E-S} that
in this case ${\left\vert\frac\beta\alpha-1\right\vert\le1}$, the
function $\dst h_0:=\left(\frac{1-z}{1+z}\right)^{\beta/\alpha}$
is univalent on $\Delta$, and its image $W=h_0(\Delta)$ is a
spiral wedge with vertex at the origin. So, by Lemma~\ref{5},
there is a point $\eta\in\pl\Delta$ such that $h(\eta)=\infty$ and
$Q_h(\eta)$ exists finitely with $\arg Q_h(\eta)=\arg\beta$ and
$\left\vert\frac\beta\alpha\right\vert\leq\left\vert
Q_h(\eta)\right\vert$.

Finally, we note that $\varphi(z)=h^{-1}(h_0(z))$ is, in fact, a
univalent function on $\Delta$. Now, applying Lemma~\ref{2} with
$\mu=\beta$, we complete the proof of the implication
(iii)$\Longrightarrow $(i) of Theorem~1, as well as assertion (c)
of Theorem~2.

Theorems 1 and 2 are proved. \epr

\noindent{\bf Proof Theorem 3.} We already know by (\ref{15-b})
and (\ref{aux5}) that ${\varphi=h^{-1}\circ h_0}$, where $h$ is
the spirallike (starlike) function associated to $f$ and\linebreak
${h_0(z))=e^{i\theta}\left(\frac{1-z}{1+z}\right)^{\beta/\alpha}}$
with $\beta=f'(\tau),\ \Re\beta>0$ and $\alpha\ge-\gamma$. So, by
Definition~\ref{def4},
\[
h_0(\Delta)=W_{\frac\beta\alpha,\theta}\subset h((\Delta).
\]
Thus
$\Omega=\varphi(\Delta)=h^{-1}\left(W_{\frac\beta\alpha,\theta}\right)$
is maximal if and only if the spiral wedge
$W_{\frac\beta\alpha,\theta}$ is maximal. In turn, by
Lemma~\ref{5}, this wedge $W_{\frac\beta\alpha,\theta}$ is maximal
if and only if $\frac\beta\alpha=-\nu$. Comparing this fact with
(\ref{aux6}), we obtain the equivalence of assertions (i) and (ii)
of the theorem.

We prove the equivalence of assertions (ii) and (iii) for the case
where $\tau=1$. Namely, let $f_1\in\G^+[1,\eta]$ with
$f_1'(1)=\beta_1>0$ and $f_1'(\eta)=\gamma_1<0$. Let $\psi$ be a
univalent solution of equation (\ref{f**}), i.e.,
\be\label{aux7}
\alpha \psi'(z)(z^{2}-1)=2f(\psi(z))
\ee
for some $\alpha\ge-\gamma_1$, normalized by $\psi(1)=1,\
\psi(-1)=\eta$.

Substituting in formula (\ref{15-b}) the explicit form of $h_0$
(see (\ref{aux5})) and the integral representation~(\ref{int})
with $\tau=1$ for the spirallike function $h$ and taking into
account that $Q_h(\eta)=\nu=\frac{\beta_1}{\gamma_1}$ (cf.
(\ref{aux6})), we get
\bep
&&(\psi(z)-1)(1-\psi(z)\bar\eta)^{\beta_1/\gamma_1}\cdot\\ &&\cdot
\exp\left[-\left(2+\beta_1/\gamma_1\right)\int_{\pl\Delta}\log(1-\psi(z)\bar\zeta)d\sigma(\zeta)
\right] =C_1\left(\frac{1-z}{1+z}\right)^{\beta_1/\alpha}
\eep
or
\bep
&&\frac{\psi(z)-\eta}{z+1}=
(z+1)^{-1-\gamma_1/\alpha}\cdot\,\frac{C_1(1-z)^{\gamma_1/\alpha}}{(1-\psi(z))^{\gamma_1/\beta_1}}\cdot\\
&&\cdot\exp\left[\frac{2\gamma_1+\beta_1}{\beta_1}\,
\int_{\pl\Delta}\log(1-\psi(z)\bar\zeta)d\sigma(\zeta) \right].
\eep
Note that one can choose an analytic branch of the multivalued
function
$C_1\dst\,\frac{(1-z)^{\gamma_1/\alpha}}{(1-\psi(z))^{\gamma_1/\beta_1}}$.
We denote this branch by $\chi(z)$. It is a continuous function
which does not vanish at the point $z=-1$. Hence its argument is a
well-defined continuous function at this point. Thus
\bep
&& \arg\frac{\psi(z)-\eta}{z+1} =\\ && =\arg\left(
(z+1)^{-1-\gamma_1/\alpha} \right)+\arg\chi(z)
+\frac{2\gamma_1+\beta_1}{\beta_1}
\int_{\pl\Delta}\arg(1-\psi(z)\bar\zeta)d\sigma(\zeta) .
\eep

Exactly as in the proof of Lemma~\ref{4}, we conclude that the
limit of the last summand exists finitely. Therefore, the function
$\psi$ is isogonal if and only if the limit
\[
\lim_{z\to-1}\arg\left((z+1)^{-1-\gamma_1/\alpha}\right)
\]
exists. Obviously, this happens if and only if the exponent
vanishes, i.e., $\alpha=-\gamma_1$.

Now let $\tau\in\overline\Delta$ be arbitrary, and let $f\in\gtep$
with $f'(\tau)=\beta,\ \Re\beta>0$, and $f'(\eta)=\gamma<0$. Let
$\varphi$ be a univalent solution of equation (\ref{f**}) for some
$\alpha\ge-\gamma$, normalized by $\varphi(1)=\tau,\
\varphi(-1)=\eta$. Denote by $h$ the spirallike function
associated to $f$, that is, $h$ satisfies equation (\ref{stareq})
with $\mu=\beta$. As above, let $h_0$ be the function which maps
the disk $\Delta$ onto a spiral wedge, namely,
$h_0(z)=e^{i\theta}\left(\frac{1-z}{1+z}\right)^{\beta/\alpha}$,
such that $\varphi=h^{-1}\circ h_0$.

Repeating the constructions in the proof of Lemma~\ref{5}, we find
a conformal mapping $\Phi$ of $\Delta$ such that $\Phi(1)=\tau,\
\Phi(-1)=\eta$, and $h_1=h\circ\Phi$ is a spirallike function with
respect to a boundary point. Note here that the domain
$D=\Phi(\Delta)$ has a corner of opening $\pi$ at the point $\eta$
because $\Phi$ maps a circular arc containing $z=-1$ onto a
circular arc which contains $z=\eta$. By Theorem~3.7 of
\cite{PC1}, the limit
$\lim\limits_{z\to-1}\arg\frac{\Phi(z)-\eta}{z+1}$ exists. Hence
$\Phi$ is isogonal at the point $-1$. Moreover, by
Proposition~4.11 of \cite{PC1}, the function $\Phi$ satisfies the
Visser--Ostrowski condition
\be\label{aux8}
\angle\lim_{z\to-1}\frac{\Phi(z)-\eta}{z+1}=1.
\ee

Now write
\be\label{aux10}
\varphi=h^{-1}\circ h_0=\Phi\circ\left( h_1^{-1}\circ
h_0\right)=\Phi\circ\psi,
\ee
where $\psi=h_1^{-1}\circ h_0$. One sees that
\[
\psi(-1):=\lim_{s\to-1^+}\psi(s)=\lim_{s\to-1^+}h_1^{-1}\left(h_0(s)\right)=-1
\]
and
\[
\psi(1):=\lim_{s\to1^-}\psi(s)=\lim_{s\to1^-}h_1^{-1}\left(h_0(s)\right)=1.
\]
Using this notation, we have
\be\label{aux12}
\arg\frac{\varphi(z)-\eta}{z+1}=
\arg\frac{\Phi(\psi(z))-\eta}{\psi(z)+1}+\arg\frac{\psi(z)+1}{z+1}\,.
\ee
Thus (\ref{aux8}) and (\ref{aux12}) imply that $\varphi$ is
isogonal at the point $\eta$ if and only if  $\psi$ is isogonal at
the point $z=-1$.

Now we check that function $\psi$ satisfies equation (\ref{aux7}).
We have seen already in the proof of Lemma~\ref{5} that $\beta
h_1(z)=h_1'(z)f_1(z)$, where ${f_1\in\G[1,-1]}$ is defined by
$f_1(z)=\frac{f(\Phi(z))}{\Phi'(z)}$. Using (\ref{aux8}), we get
\[
f_1'(-1)=\angle\lim_{z\to-1}\frac{f_1(z)}{z+1}=
\angle\lim_{z\to-1}\frac{f(\Phi(z))}{\Phi(z)-\eta}\cdot\frac{\Phi(z)-\eta}{(z+1)\Phi'(z)}=\gamma.
\]
Furthermore,
\bep
h_0(z)=h_1(\psi(z))=\frac1\beta\,h_1'(\psi(z))f_1(\psi(z))\\
=\frac1\beta\,\frac{\left(h_1(\psi(z))\right)'}{\psi'(z)}\,f_1(\psi(z))=
\frac{h_0'(z)}{\beta\psi'(z)}\,f_1(\psi(z)).
\eep
Substituting
$h_0(z)=e^{i\theta}\left(\frac{1-z}{1+z}\right)^{\beta/\alpha}$ in
the last equality and differentiating, we see that equation
(\ref{aux7}) holds. But we have already shown that in that case
$\psi$ (hence, $\varphi$) is isogonal if and only if
$\alpha=-f_1'(-1)=-\gamma$. This completes the proof. \epr

\noindent{\bf Proof of Theorem 4.} Assertions (i) and (ii) of the
theorem are direct consequences of Lemma 6. To prove assertion
(iii), we first note that the inclusion $\tau\in\cap_k\pl\Omega_k$
follows by assertion (a) of Theorem~2.

Also observe that for each pair $k_1$ and $k_2$ such that
$\eta_{k_1}\neq\eta_{k_2}$, the set
$\Omega_{k_1,k_2}=\Omega_{k_1}\cap\Omega_{k_2}$ is empty. Indeed,
otherwise $\Omega_{k_1,k_2}$ is a FID for $\Ss$. Hence, it must
contain a point $\eta\in\pl\Omega_{k_1,k_2}\cap\pl\Delta$ such
that
\[
\eta=\angle\lim_{t\rightarrow -\infty }F_{t}(z)
\]
whenever $z\in\Omega_{k_1,k_2}.$ Hence we should have a
contradiction $\eta =\eta_{k_1}=\eta_{k_2}.$

Let us suppose now that for a pair $k_1$ and $k_2$ there is a
point $z_0\neq\tau$, $z_0\in\Delta$, such that
$z_0\in\pl\Omega_{k_1}\cap\pl\Omega_{k_2}.$ Then the whole curve
\[
l=\{z\in\Delta:\ z=F_t(z_0),\ t\geq 0\}
\]
ending at $\tau$ must belong to both $\overline{\Omega}_{k_1}$ and
$\overline{\Omega}_{k_2}$, hence to
$l\subset\pl\Omega_{k_1}\cap\pl\Omega_{k_2}$, since
$\Omega_{k_1}\cap\Omega_{k_2}=\emptyset.$

Finally, we have that $f\in\Hol(\Delta,\C)$ is locally
Lipshitzian. Therefore, if $\zeta\in\Delta$ is an interior end
point of $l,\ \zeta\neq\tau,$ then there is $\delta >0$ such that
the Cauchy problem (\ref{2}) has a solution
$u(t,\zeta)\left(=F_t(\zeta)\right)$ for all
$t\in[-\delta,\infty)$; and the curve $l_1=\{z\in\Delta:\
z=u(t,\zeta),\ t\in[-\delta,\infty)\}$ also belongs to
$\pl\Omega_{k_1}\cap\pl\Omega_{k_2}$. But $l_1$ properly contains
$l$, which is impossible. So $\zeta$ must belong to $\pl\Delta.$
The corollary is proved. \epr

\noindent{\bf Remark 6.} The complete solution to the problem of
finding FID's requires the treating the case in which
$\tau\in\pl\Delta$ and $f'(\tau)=0$. In this case, the semiflow
$\Ss=\{F_{t}\}_{t\geq 0}$ generated by $f$ consists of
self-mappings of $\Delta$ of parabolic type. This delicate
question is equivalent to the following problem. Associate with
$f$ a univalent function $h\in\Hol(\Delta,\C)$ which is a solution
of Abel's functional equation
\be\label{*}
h(F_{t}(z))=h(z)+Kt,\quad t\geq 0,
\ee
for some $K\in\C$ which does not depend on $t\geq 0.$ Under what
conditions does the image $h(\Delta)$ contain a strip $W$ such
that equation (\ref{*}) holds for all $t\in\R$, whenever
$z\in\Omega =h^{-1}(W)$? We hope to consider this problem
elsewhere.

\vspace{3mm}

\noindent{\bf Appendix.} Quoting T. Harris \cite{HT}, we note that
``a classical problem of analysis is a problem that has interested
mathematicians since the time of Abel: how to define the $n$-th
iterate of a function when $n$ is not an integer."

In other words, the question is given a function $F$, to find a
family of functions $\left\{F_{t}\right\}_{t\geq 0},$ with
$F_{1}=F$ satisfying the semigroup (group) property for all $t\geq
0$ ( respectively, $t\in\R$). This problem is called the embedding
problem into a continuous semiflow (respectively, flow).

The possibility of such an embedding is important, in particular,
in problems of conformal mapping and in the study of Markov
branching processes with continuous time (whose first general
formulation appears to have been given by Kolmogorov (1947)).

When $F$ is a holomorphic function, K{\oe}nigs (1884) showed how
the problem may be solved locally near a fixed point $z_{0}$ such
that ${0<|F'(z_{0})| <1}$.

The limit
\[
\lim_{n\to\infty}\frac{F^{n}(z_{0})-z_{0}}{(F'(z_{0}))^{n}}=h(z)
\]
can be shown to exist for $z$ near $z_{0}$ and to satisfy
Schr\"oder's functional equation
\begin{equation}\label{N}
h(F(z))=F^{\prime }(z_{0})h(z),
\end{equation}
whence
\[
F(z)=h^{-1}[F'(z_{0})h(z)].
\]

The latter expression then serves as a definition of $F_{t}$ when
$t$ is not necessarily an integer:
\[
F_{t}=h^{-1}[(F^{\prime }(z_{0}))^{t}h(z)].
\]
Consequently, if $F\in\Hol(\Delta )$ is a self-mapping of the unit
disk $\Delta$ and $z_{0}\in \Delta $, then
$\Ss=\{F_{t}\}_{t\geq0}$ is {\bf globally} well-defined on
$\Delta$ if and only if $h$ is a $\mu$-spirallike function on
$\Delta$ with $\arg\mu=\arg(-\log F'(z_{0})).$

Following the work of I. N. Baker \cite{BIN}, S. Karlin and J.
McGregor \cite{K-MG} considered the local embedding problem of
holomorphic functions with two fixed points into a continuous
group. In particular, they studied a class $\Ll$ of functions
holomorphic in the extended complex plane $\overline\C$ except for
an at most countable closed set in $\overline\C$ and proved the
following result.

{\it Let $F$ be a function of class $\Ll$ with two fixed points
$z_{0}$ and $z_{1}$, such that the segment $[z_{0},z_{1}]$ is in
the domain of regularity of $F$ and is mapped onto itself. Assume
that $0<|F'(z_{0})| <1<|F'(z_{1})|$ and that for $z$ in the open
segment $(z_{0},z_{1}),$ $F(z)\neq z,$ $F^{\prime }(z)\neq 0.$
Then there is a continuous one-parameter group
$\{F_{t}\}_{t\in\R}$ of functions with common fixed points $z_{0}$
and $z_{1}$ and invariant segment $[z_{0},z_{1}]$ such that
$F_{1}(z)=F(z)$ if and only if $F(z)$ is a linear fractional
transformation on $\mathbb{C}.$}

First we note that the condition that $F$ map $[z_{0},z_{1}]$ into
itself implies that $F^{\prime }(z)$ is real on this segment.

Suppose now that $F$ is linear fractional, $F(z)\not\equiv z$, and
let $z_{0}$ and $z_{1}$ be its finite fixed points, $z_0\not=z_1$.
The following simple assertion can be obtained by using the linear
model of mappings having two fixed points $0$ and $\infty$ and
applying the Julia--Carath\'eodory theorem.

\begin{lemma}
The following are equivalent.

(i) There is an open disk $D$ such that either $z_{0}\in\pl D$ and
$z_1\not\in\overline{D}$, or $z_0\in D$ and $z_1\in\pl D$, which
is $F$-invariant.

(ii) Each open disk $D$ such that $z_{0}\in\overline{D}$ and
$z_{1}\notin D$ is $F$-invariant.

(iii) The segment $[z_{0},z_{1}]$ is $F$-invariant and
$|F'(z_0)|\le1$.

(iv) If $a=F'(z_{0})$ then $0<a<1.$
\end{lemma}

Since Schr\"oder's equation (\ref{N}) with linear-fractional $F$
has a linear-fractional solution $h$, we have that $h$ is
starlike; hence $F$ can be embedded into a one-parameter semigroup
$\{F_{t}\}_{t\in\R}$ on each disk $D$ containing $z_{0}$ and such
that $z_{1}\notin D$. This disk is $F_{t}$-invariant for all
$t\geq 0$.

In turn, for the embedding property into a continuous group, we
obtain the following assertion by using our Theorems~1 and~2 and
Theorem~1 in \cite{K-MG}.

\begin{corol}
Let $F$ be a function of class $\Ll$ with two different fixed
points $z_0$ and $z_1$. Assume that $0<|F'(z_0)| <1<|F'(z_1)|$,
and that for $z$ in the open segment $(z_0,z_1),$ $F(z)\neq z,$
$F'(z)\neq 0$. The following assertions are equivalent.

(i) For each open disk $D$ such that $z_0\in\overline{D}$ and
$z_1\notin D,$ there is a semiflow $\Ss=\left\{F_t\right\}_{t\geq
0}$ with $F_1=F$ such that $\Ss\subset\Hol(D)$.

(ii) For each domain $\Omega$ bounded by two circles passing
through $z_0$ and $z_1$, there is a one-parameter flow
$\Ss=\{F_t\}_{t\in\mathbb{R}}$ such that $\Ss\subset\Aut(\Omega)$
and $F=F_1.$

(iii) The function $F$ is linear fractional with
$0<F^{\prime}(z_{0})<1.$

Consequently, in this case, for any disk $D$ such that
$z_0\in\overline D$ and $z_1\in\pl D$, the maximal (backward)
flow-invariant domain is the disk\ $\Omega\subset D$ whose
boundary passes through $z_0$ and is internally tangent to $\pl D$
at $z_1$.
\end{corol}


\end{document}